\documentclass[a4paper,11pt,leqno]{article}
\usepackage{amsmath}
\usepackage{amsfonts}
\usepackage{amsthm}
\usepackage{amssymb}
\usepackage{amscd}

\title{General elements in anticanonical systems of
threefolds with divisorial contractions
and applications to classification}
\author{Masayuki Kawakita}
\date{}

\theoremstyle{definition}
\newtheorem{Definition}{Definition}[section]
\newtheorem{Example}[Definition]{Example}
\newtheorem{Text}[Definition]{}

\theoremstyle{plain}
\newtheorem{Theorem}[Definition]{Theorem}
\newtheorem{Proposition}[Definition]{Proposition}
\newtheorem{Lemma}[Definition]{Lemma}
\newtheorem{Corollary}[Definition]{Corollary}

\numberwithin{equation}{Definition}

\theoremstyle{definition}
\newtheorem{sText}[equation]{}
\newtheorem{Remark}[equation]{Remark}

\theoremstyle{plain}

\newtheorem{sLemma}[equation]{Lemma}

\newcommand{\bC}{{\mathbb C}}
\newcommand{\bP}{{\mathbb P}}
\newcommand{\bQ}{{\mathbb Q}}
\newcommand{\bZ}{{\mathbb Z}}

\newcommand{\cI}{{\mathcal I}}

\newcommand{\cL}{{\mathcal L}}

\newcommand{\cO}{{\mathcal O}}
\newcommand{\cQ}{{\mathcal Q}}
\newcommand{\cR}{{\mathcal R}}

\newcommand{\fm}{{\mathfrak m}}

\newcommand{\divis}{\mathrm{div}}
\newcommand{\Ext}{\mathrm{Ext}}
\newcommand{\Hom}{\mathrm{Hom}}

\newcommand{\red}{\mathrm{red}}
\newcommand{\Sing}{\mathrm{Sing}}

\newcommand{\wt}{\mathrm{wt}}

\begin{document}
\maketitle

\begin{abstract}
We treat threefolds with divisorial contractions whose
exceptional divisors contract to compound Du Val points.
We prove that general elements in their anticanonical systems
around the exceptional divisors have at worst Du Val singularities.
As applications to classification, we describe divisorial contractions
to compound $A_n$ points, and moreover we conclude that discrepancies of
divisorial contractions to compound $D_n$ or $E_n$ points
are at most four.
\end{abstract}

\section{Introduction}\label{sec:intro}
This paper aims at completion of the explicit study of three dimensional
divisorial contractions whose exceptional divisors contract to
compound Du Val points after \cite{Ka01} and \cite{Ka01p}.

It is M.\@ Reid who has pointed out that general elements
in the anticanonical systems of threefolds have at worst
Du Val singularities under appropriate situations
arising from contractions of extremal faces, and it has pervaded us
as the general elephant conjecture.
\cite{Sh79}, \cite{Re83p} and \cite{Ta00} sustain it by affirmative
answers for Fano threefolds with singularities and
this approach settles the existence problem of three dimensional flips
\cite{Mo88}\cite{KoMo92}.

Returning to our case of divisorial contractions,
let $f \colon (Y \supset E) \to (X \ni P)$ be a germ of
a three dimensional divisorial contraction
whose exceptional divisor $E$ contracts to a point $P$.
His conjecture claims that
general elements in the anticanonical system $\left| -K_Y \right|$ of $Y$
around $E$ have Du Val singularities only.
This statement is analogue to the flipping case and
we can start at a similar stage
once we take the intersection $C$ of $E$ and
the strict transform of a general hyperplane section on $X$,
because the first cohomology of $C$ vanishes.

There lies, however, a crucial difference in $C$
between our case and the flipping case.
For any flipping curve in a threefold the first cohomology of
the structure sheaf of any closed subscheme supported on this curve
always vanishes, but it does not hold for our curve $C$.
Nevertheless provided that $P$ is Gorenstein, or equivalently, compound Du Val,
numerical information on $f$ in \cite{Ka01} benefits
analysis of the local structure of $C \subset Y$
through delicate attention to the way of local embedding of $Y$
into the tangent space at each non-Gorenstein point.
It stores the behaviour of global sections
in $\left| -K_Y \right|$, and finally concludes the existence of
Du Val sections (\ref{thm:gen.element}).

\begin{Theorem}\label{thm:introDuval}
Let $f\colon (Y \supset E) \to (X \ni P)$ be a germ of a
divisorial contraction whose exceptional divisor $E$ contracts
to a cDV point $P$.
Then a general element in $\left| -K_Y \right|$
has at worst Du Val singularities.
\end{Theorem}

Let $S$ be a surface on $Y$ defined by a general element
in $\left| -K_Y \right|$ and let $S_X$ be its strict transform on $X$.
Our main theorem guarantees that $S$ and $S_X$ have Du Val singularities only
and that the induced morphism $S \to S_X$ factors the minimal resolution of
$S_X$. We show it together with providing information
on the partial resolution $S \to S_X$ (\ref{thm:detail}).
In some special cases, to be precise in types O and I in
(\ref{thm:num.classn}), we see immediately that a general hyperplane section
$S_X$ on $X$ gives a Du Val section $S$.
Considerably delicate analysis is required in the remaining cases,
but thanks to this analysis we can obtain
the strong version of the general elephant conjecture, which asserts that
the type of a Du Val singularity $Q \in S$ for any point $Q$ is
nothing but that of a general Du Val section of a germ $Q \in Y$.
In most cases it is proved by the direct search for $S$, whereas
in some exceptional cases where $-K_Y$ is linearly equivalent to
the sum of $E$ and a Cartier divisor $L$ we obtain a desired $S$ by
showing that $E$ has a Du Val singularity of the required type at any $Q \in S$
and that $L$ has no base points. We add that complete description of
the partial resolution $S \to S_X$ is also obtained. For instance,
its exceptional locus is irreducible in nearly every case.

We apply this theorem to classification of
divisorial contractions naturally.
We describe these contractions to c$A_n$ points (\ref{thm:cAn}),
following smooth and c$A_1$ case \cite{Ka01}\cite{Ka01p}.

\begin{Theorem}
Assume that $P$ is c$A_n$ \textup{(}$n \ge 2$\textup{)}.
Then one of the following holds.
\begin{enumerate}
\item Under a suitable identification
$P \in X \cong o \in (x_1x_2+g(x_3, x_4)=0)
\subset \bC^4=x_1x_2x_3x_4\textrm{-}\mathrm{space}$,
$f$ is the weighted blowup with weights
$\wt(x_1, x_2, x_3, x_4)=(r_1,r_2,a,1)$,
where $a$ divides $r_1+r_2$ and is coprime to $r_1$ and $r_2$,
the weighted order of $g$ with weights $\wt(x_3, x_4)=(a,1)$
is $r_1+r_2$, and the coefficient of $x_3^{(r_1+r_2)/a}$ in $g$ is not zero.
Moreover any such $f$ is a divisorial contraction.
\item $P$ is a c$A_2$ point isomorphic to
$o \in (x_1x_2+x_3^3+g_{\ge4}(x_3, x_4)=0)
\subset \bC^4=x_1x_2x_3x_4\textrm{-}\mathrm{space}$,
where the total order of $g_{\ge4}$ with respect to $x_3, x_4$
is greater than or equal to $4$,
$Y$ has exactly one non-Gorenstein point $Q$, which is
isomorphic to $o \in (y_1^2+y_2^2+y_3^2+y_4^3=0)$ in the quotient space of
$\bC^4=y_1y_2y_3y_4\textrm{-}\mathrm{space}$ divided by $\bZ/(4)$ with
weights $\wt(y_1, y_2, y_3, y_4)=(1, 3, 3, 2)$, and $K_Y=f^*K_X+3E$.
Moreover there exists such an example.
\end{enumerate}
\end{Theorem}

It is anticipated that the number of divisorial contractions over $P$
decreases as the singularity of $P$ becomes worse on account of restriction
to choice of coordinates at $P$. Nevertheless we might not
expect that it makes the perfect explicit study of them simple because
the defining equation of $X$ at $P$ becomes much more complicated.
By this reason, instead of complete description,
we restrict the possibility of divisorial contractions
in the remaining c$D_n$ and c$E_n$ cases
by giving an upper bound of discrepancies (\ref{cor:cDcE}),
equipped with some typical examples.

\begin{Theorem}
Assume that $P$ is c$D_n$ or c$E_n$, and let $K_Y=f^*K_X+aE$. Then $a \le 4$.
\end{Theorem}

We now have sufficient tools toward classification of divisorial contractions
even in c$D_n$ and c$E_n$ cases once $P \in X$ is given explicitly.
Basically what we should do is to compare discrepancies as in \cite{Km96}
or in (\ref{lem:gen.mtd}), with enormous information on multiplicities
along $E$ of surfaces on $X$ with special directions provided by
the singular Riemann-Roch technique and the general elephant theorem,
occasionally applying Shokurov's connectedness lemma as in
\cite[Theorem 3.10]{Co00}, \cite[Theorem 3.6]{CM00}.

This paper is constructed as follows.
In Section \ref{sec:preliminaries} we state our theorems precisely.
Section \ref{sec:numerical} is devoted to preparing basic numerical techniques.
After recalling results for divisorial contractions
to compound Du Val points in \cite{Ka01}, we introduce
our fundamental setup toward analysis of $C \subset Y$ following \cite{Mo88}.
Quite delicate local investigation of $C \subset Y$ is
presented in Section \ref{sec:local}.
Using it in Section \ref{sec:existence} we prove
our main theorem, the existence of Du Val sections.
In Section \ref{sec:possible} we restrict possible divisorial contractions
focusing on the types of singularities on Du Val sections.
It provides an upper bound of discrepancies in c$D_n$ and c$E_n$ cases.
Finally in Section \ref{sec:cAn} we give explicit description
of divisorial contractions in c$A_n$ case.

I would like to thank Professor Yujiro Kawamata
for his stimulating encouragement.
I am grateful to Professor Alessio Corti and Professor Miles Reid
for motivating me in this subject.
I have written this paper during my visit to the University of Cambridge,
where I receive warm hospitality of Professor Corti.
Financial support has been provided by
the Japan Society for the Promotion of Science.

\section{Preliminaries and statements}\label{sec:preliminaries}
We work over the complex number field $\bC$.
We fundamentally do in the analytic category but sometimes enter
the algebraic category through algebraisation theorems
of M.\@ Artin \cite{Ar68}\cite{Ar69}.
First we define a divisorial contraction in a general sense.

\begin{Definition}
Let $f \colon Y \to X$ be a morphism with connected fibres between
normal varieties with at worst terminal singularities.
We call $f$ a \textit{divisorial contraction}
if the exceptional locus of $f$ is a prime divisor and $-K_Y$ is $f$-ample.
\end{Definition}

\begin{Remark}
Of course this definition applies to divisorial contractions
emerging in the usual minimal model program \cite{KMM87}.
\end{Remark}

We recall the classification of three dimensional terminal singularities.

\begin{Definition}
Let $P \in X$ be a germ of a three dimensional variety.
We call $P$ a \textit{cDV} (\textit{compound Du Val}) \textit{point}
if a general hyperplane section has at worst a Du Val singularity at $P$.
The singularity $P$ is said to be
\textit{c$A_n$, c$D_n$, c$E_n$} (\textit{compound $A_n, D_n, E_n$})
according to the type of the Du Val singularity on
a general hyperplane section.
\end{Definition}

\begin{Remark}
We say that a smooth point on a surface or a threefold is
$A_0$, c$A_0$ respectively for convenience.
\end{Remark}

\begin{Theorem}[{\cite[Theorem1.1]{Re83}}]
Let $P \in X$ be a germ of a three dimensional variety.
Then $P$ is a Gorenstein terminal singularity if and only if
$P$ is an isolated cDV point.
\end{Theorem}

\begin{Text}
Consider a non-Gorenstein terminal singularity $P \in X$.
Let $r$ be the local Gorenstein index of $P \in X$,
that is, the smallest positive integer such that $rK_X$ is Cartier at $P$.
Take the index one cover $\pi \colon (X^\sharp \ni P^\sharp) \to (X \ni P)$,
which is a cyclic $\mu_r$-cover.
Fix a character generating $\Hom(\mu_r, \bC^{\times})=\bZ/(r)$ and
define the weight modulo $r$ for any semi-invariant function on $X^\sharp$
with respect to this character.
\end{Text}

\begin{Theorem}[{\cite{Mo85}}]\label{thm:non.Gor}
There exists a $\mu_r$-equivariant identification
\begin{align*}
P^\sharp \in X^\sharp \cong o \in (\phi=0) \subset \bC^4=
x_1x_2x_3x_4\textrm{-}\mathrm{space},
\end{align*}
where $x_1, x_2, x_3, x_4$ and $\phi$ are $\mu_r$-semi-invariant
and $\phi=x_4$ if $P^\sharp \in X^\sharp$ is a smooth point.
The weights of $x_1, x_2, x_3, x_4$ and $\phi$ satisfy one of the following.
\begin{enumerate}
\item $\wt(x_1, x_2, x_3, x_4; \phi)=(1, -1, b, 0; 0)$, where
      $b$ is coprime to $r$.\label{itm:ordin}
\item $r=4$ and
$\wt(x_1, x_2, x_3, x_4; \phi)=(1, 3, 3, 2; 2)$.\label{itm:excep}
\end{enumerate}
\end{Theorem}

\begin{Remark}
\cite{Mo85} gives more precise description of $\phi$ and weights,
and the classification is completed by \cite[Theorem 6.4]{KoSB88}.
\end{Remark}

\begin{Remark}\label{rem:basket}
Any three dimensional non-Gorenstein terminal singularity $P \in X$ has
a small deformation to a basket of terminal quotient singularities $P_i$.
$P_i$ are called fictitious singularities in the sense of M.\@ Reid
\cite{Re87}.
If $P$ is of type (\ref{itm:ordin}) in (\ref{thm:non.Gor}) then
any local index at $P_i$ equals that at $P$.
If $P$ is of type (\ref{itm:excep}) then
the local index at one of $P_i$ equals $2$ and that at each of
the rest equals $4$.
\end{Remark}

Now it is the time when we state theorems.
To begin with we give a numerical classification of $f$
which is a slightly revised version of \cite[Theorem 4.5]{Ka01}
quoted in (\ref{thm:num.classn}), due to (\ref{thm:O.I.elephant}).

\begin{Theorem}\label{thm:re.num.classn}
Let $f \colon (Y \supset E) \to (X \ni P)$ be a germ of
a three dimensional divisorial contraction
whose exceptional divisor $E$ contracts to a cDV point $P$.
Let $K_Y=f^*K_X+aE$,
let $I = \{Q: \textrm{type} \ \frac{1}{r_Q}(1, -1, \overline{av_Q})\}$
\textup{(}$v_Q \le r_Q / 2$\textup{)}
be the set of fictitious singularities
from non-Gorenstein singularities on $Y$,
let $J=\{(r_Q, v_Q)\}_{Q \in I}$,
and let $1+d(-1)=\dim \cO_X/f_*\cO_Y(-2E)$.
Precise notation is provided in \textup{(\ref{txt:notn})}.
Then $f$ is exactly one of the following types.
\begin{center}
\begin{tabular}{l|c|l|l}
\hline
\multicolumn{1}{c|}{type} & $1+d(-1)$ &
\multicolumn{1}{c|}{$J$}  & \multicolumn{1}{c}{$a$} \\
\hline
O    & $\ge 2$ &                          & $1$ \\
I    & $1$     & $\{(7, 3)\}$ or $\{(3, 1), (5, 2)\}$ & $2$ \\
IIa  & $2$     & $\{(r, 2)\}$             & $2$ or $4$ \\
IIb  & $2$     & $\{(r_1, 1), (r_2, 1)\}$ & $(r_1+r_2)/r_1r_2E^3 \ge2$ \\
III  & $3$     & $\{(r, 1)\}$             & $(1+r)/rE^3 \ge2$ \\
IV   & $4$     & $\emptyset$              & $2$ \\
\hline
\end{tabular}
\end{center}
\end{Theorem}

\begin{Remark}\label{rem:divide}
For convenience we divide type IIb into two types according to
the number of non-Gorenstein points on $Y$.
We say that $f$ is of type ${\textrm{IIb}}^{\vee}$
or ${\textrm{IIb}}^{\vee \vee}$ if $f$ is of type II
with one, two non-Gorenstein points on $Y$ respectively.
\end{Remark}

Our main theorem is on general elements in the anticanonical system
$\left| -K_Y \right|$
of $Y$, dubbed general elephants by M.\@ Reid \cite{Re87}.

\begin{Theorem}\label{thm:gen.element}
Let $f\colon (Y \supset E) \to (X \ni P)$ be a germ of a 
divisorial contraction in \textup{(\ref{thm:re.num.classn})}.
Then a general element in $\left| -K_Y \right|$
has at worst Du Val singularities.
\end{Theorem}

Let $S$ be a surface on $Y$ defined by a general element
in $\left| -K_Y \right|$ and let $S_X$ be its strict transform on $X$.
$S$ and $S_X$ have at worst Du Val singularities and
the induced morphism $S \to S_X$ factors the minimal resolution of $S_X$
by (\ref{thm:gen.element}).
We obtain the main theorem together with information
on the partial resolution $S \to S_X$.
(\ref{thm:gen.element}) and (\ref{thm:detail}) are induced
by (\ref{thm:O.I.elephant}), (\ref{thm:el.irred}), (\ref{thm:II.III.elephant}),
(\ref{thm:typeI}), (\ref{thm:typeIIandIII}) with (\ref{rem:typeIV})
and (\ref{rem:typeIIandIII}).

\begin{Theorem}\label{thm:detail}
\begin{enumerate}
\item Assume that $f$ is of type O or I.
Then $S_X$ defines a general element in $\left| -K_X \right|$.
Moreover if $f$ is of type I
then $P$ is c$E_7$ or c$E_8$.
\item Assume that $f$ is of type II, III or IV.
Then possible types of Du Val singularities on $S$ and $S_X$
are as follows.
\begin{center}
\begin{tabular}{l|l|l}
\hline
\multicolumn{1}{c|}{type of $f$}   &
\multicolumn{1}{c|}{type of $S_X$} & \multicolumn{1}{c}{type of $S$} \\
\hline
IIa  & $D_r$ or $D_{r+1}$                       & $A_{r-1}$ \\
${\textrm{IIb}}^{\vee}$, $J=\{(r, 1), (r, 1)\}$
& $D_{2r}$ or $D_{2r+1}$ & $A_{2r-1}$ \\
${\textrm{IIb}}^{\vee}$, $J=\{(3, 1), (3, 1)\}$ & $E_7$    & $E_6$ \\
${\textrm{IIb}}^{\vee}$, $J=\{(2, 1), (4, 1)\}$ & $E_6$    & $D_5$ \\
${\textrm{IIb}}^{\vee \vee}$ & $A_{r_1+r_2-1}$  & $A_{r_1-1}$ and $A_{r_2-1}$\\
III  & $A_r$                           & $A_{r-1}$ \\
IV   & smooth                          & smooth \\
\hline
\end{tabular}
\end{center}
The type of a Du Val singularity $Q \in S$ for any point $Q$ is
that of a general Du Val section of a germ $Q \in Y$.
Furthermore the exceptional locus of the partial resolution $S \to S_X$
is irreducible unless the type of $S_X$ is that after
the word ``or" in the above table.
\end{enumerate}
\end{Theorem}

\begin{Remark}
We have examples of type I with $J=\{(3, 1), (5, 2)\}$ in (\ref{exl:I})
but I do not know whether type I with $J=\{(7, 3)\}$ happens or not.
If $f$ is of type I with $J=\{(7, 3)\}$ then $P$ is c$E_7$ (\ref{rem:typeI}).
\end{Remark}

We bound the discrepancy $a$ in the case where $P$ is c$D_n$ or c$E_n$
as a corollary.

\begin{Corollary}\label{cor:cDcE}
Assume that $P$ is c$D_n$ or c$E_n$.
Then $f$ is of type O, I, IIa or
${\textrm{IIb}}^{\vee}$ in \textup{(\ref{thm:re.num.classn})}
and the discrepancy $a \le 4$.
\end{Corollary}

\begin{Remark}
We have examples of type O, I, ${\textrm{IIb}}^{\vee}$,
or with $a=1, 2, 3$ in
(\ref{exl:3sing}), (\ref{exl:I}), (\ref{exl:rr}) and (\ref{exl:24}).
However I do not know whether type IIa happens or not.
If $a=4$ in c$D_n$ or c$E_n$ case then $r=5$ and
$P$ has to be c$D_4$, c$D_5$ or c$D_6$ (\ref{thm:localIIa}(\ref{itm:th4ir})).
\end{Remark}

Divisorial contractions to smooth or c$A_1$ points
have been completely classified
in \cite{Ka01} and \cite{Ka01p}.

\begin{Theorem}[{\cite{Ka01}}]
Assume that $P$ is smooth.
Then $f$ is a weighted blowup.
More precisely, we can take local coordinates $x_1, x_2, x_3$ at $P$
and coprime positive integers $s$ and $t$,
such that $f$ is the weighted blowup of $X$ with weights
$\wt (x_1, x_2, x_3) = (1, s, t)$.
Moreover any such $f$ is a divisorial contraction.
\end{Theorem}

\begin{Remark}
We may assume that $s \le t$.
The type of $f$ is ${\textrm{IIb}}^{\vee \vee}$, III, IV
if $s>1$, $s=1$ and $t>1$, $t=1$ respectively.
\end{Remark}

\begin{Theorem}[{\cite{Ka01p}}]
Assume that $P$ is c$A_1$.
Then $f$ is a weighted blowup.
More precisely, under a suitable identification
$P \in X \cong o \in (x_1x_2+x_3^2+x_4^N=0) \subset \bC^4
=x_1x_2x_3x_4\textrm{-}\mathrm{space}$,
$f$ is the weighted blowup with one of the following weights.
\begin{enumerate}
\item $\wt(x_1, x_2, x_3, x_4)=(s,2t-s,t,1)$,
where $s,t$ are coprime positive integers such that
$s \le t \le N/2$.\label{itm:A1ord}
\item $N=3$ and $\wt(x_1, x_2, x_3, x_4)=(1,5,3,2)$.\label{itm:A1exc}
\end{enumerate}
Moreover any such $f$ is a divisorial contraction.
\end{Theorem}

\begin{Remark}
In (\ref{itm:A1ord}) the type of $f$ is ${\textrm{IIb}}^{\vee \vee}$, III, O
if $s>1$, $s=1$ and $t>1$, $t=1$ respectively.
In (\ref{itm:A1exc}) the type of $f$ is IIa.
\end{Remark}

\begin{Remark}
A.\@ Corti has obtained the result in the case where $P$
is an ordinary double point earlier \cite[Theorem 3.10]{Co00}.
\end{Remark}

The case where $P$ is c$A_n$ ($n \ge 2$) is done in Section \ref{sec:cAn}
in this paper.

\begin{Theorem}\label{thm:cAn}
Assume that $P$ is c$A_n$ \textup{(}$n \ge 2$\textup{)}.
Then one of the following holds.
\begin{enumerate}
\item Under a suitable identification
$P \in X \cong o \in (x_1x_2+g(x_3, x_4)=0)
\subset \bC^4=x_1x_2x_3x_4\textrm{-}\mathrm{space}$,
$f$ is the weighted blowup with weights
$\wt(x_1, x_2, x_3, x_4)=(r_1,r_2,a,1)$,
where $a$ divides $r_1+r_2$ and is coprime to $r_1$ and $r_2$,
the weighted order of $g$ with weights $\wt(x_3, x_4)=(a,1)$
is $r_1+r_2$, and the coefficient of $x_3^{(r_1+r_2)/a}$ in $g$ is not zero.
Moreover any such $f$ is a divisorial contraction.\label{itm:cAn.ord}
\item $P$ is a c$A_2$ point isomorphic to
$o \in (x_1x_2+x_3^3+g_{\ge4}(x_3, x_4)=0)
\subset \bC^4=x_1x_2x_3x_4\textrm{-}\mathrm{space}$,
where the total order of $g_{\ge4}$ with respect to $x_3, x_4$
is greater than or equal to $4$,
$Y$ has exactly one non-Gorenstein point $Q$, which is
isomorphic to $o \in (y_1^2+y_2^2+y_3^2+y_4^3=0)$ in the quotient space of
$\bC^4=y_1y_2y_3y_4\textrm{-}\mathrm{space}$ divided by $\bZ/(4)$ with
weights $\wt(y_1, y_2, y_3, y_4)=(1, 3, 3, 2)$,
and $K_Y=f^*K_X+3E$.
Moreover there exists such an example
\textup{(\ref{exl:cA2})}.\label{itm:cA2.exc}
\end{enumerate}
\end{Theorem}

\begin{Remark}
In (\ref{itm:cAn.ord}) we may assume that $r_1 \le r_2$.
The type of $f$ is ${\textrm{IIb}}^{\vee \vee}$, III, O
if $r_1>1$ and $a>1$, $r_1=1$ and $a>1$, $a=1$ respectively.
In (\ref{itm:cA2.exc}) the type of $f$ is ${\textrm{IIb}}^{\vee}$
with $J=\{(2, 1), (4, 1)\}$.
\end{Remark}

\begin{Remark}
A.\@ Corti and M.\@ Mella have obtained the result in the case where $P$
is a c$A_2$ point isomorphic to $o \in (x_1x_2+x_3^3+x_4^3=0)
\subset \bC^4$ earlier \cite[Theorem 3.6]{CM00}.
\end{Remark}

\section{Basic numerical results}\label{sec:numerical}
\begin{Text}\label{txt:notn}
Let $f \colon (Y \supset E) \to (X \ni P)$ be a germ of
a three dimensional divisorial contraction
whose exceptional divisor $E$ contracts to a cDV point $P$.
We remark that $f$ can be always extended to a morphism between
projective varieties \cite{Ar68}\cite{Ar69}.
Let $K_Y=f^*K_X+aE$ and let $r$ be the global Gorenstein index of $Y$.

\begin{sLemma}[{\cite[Lemma 4.3]{Ka01}}]
$a$ and $r$ are coprime.
\end{sLemma}

We take an integer $e$ such that $ae \equiv 1$ modulo $r$.
Let $I = \{Q: \textrm{type} \ \frac{1}{r_Q}(1, -1, b_Q)\}$
be the set of fictitious singularities
from non-Gorenstein singularities on $Y$ (\ref{rem:basket}).
Then $(\cO_{Y_Q}(E_Q))_Q \cong (\cO_{Y_Q}(eK_{Y_Q}))_Q$,
where $(Y_Q, E_Q)$ is the deformed pair at $Q$ from $(Y, E)$.
We note that $b_Q$ is coprime to $r_Q$ and that
$e$ is also coprime to $r_Q$ because $ae \equiv 1$ modulo $r$.
Hence $v_Q= \overline{eb_Q}$ is coprime to $r_Q$.
Here $\bar{\ }$ denotes the smallest residue modulo $r_Q$,
that is, $\overline{j} = j - \lfloor \frac{j}{r_Q} \rfloor r_Q$,
where $\lfloor \ \rfloor$ denotes the round down, that is,
$\lfloor j \rfloor = \max\{k \in \bZ \mid k \le j\}$.
Replacing $b_Q$ with $r_Q - b_Q$ if necessary,
we may assume that $v_Q \le r_Q / 2$.
With this description, $r=1$ if $I$ is empty, and otherwise
$r$ is the lowest common multiple of $\{r_Q\}_{Q \in I}$.
We set $J=\{(r_Q, v_Q)\}_{Q \in I}$.
We moreover define $d(i)=\dim f_*\cO_Y(iE)/f_*\cO_Y((i-1)E)$.
We note that $d(i) = 0$ ($i \ge 1$) and $d(0) = 1$.
\end{Text}

\begin{Remark}\label{rem:a<r}
$a < \max\{r_Q\}$ unless $P$ is a smooth point
by analogue to the proof of \cite[Lemma 6.10]{Ka01p},
because there exists a valuation with centre $P$ whose
discrepancy with respect to $K_X$ is $1$ \cite{Ma96}.
\end{Remark}

The next theorem gives a rough numerical classification of $f$.

\begin{Theorem}[{\cite[Theorem 4.5]{Ka01}}]\label{thm:num.classn}
$f$ is exactly of one of the following types.
\begin{center}
\begin{tabular}{l|c|l|l}
\hline
\multicolumn{1}{c|}{type} & $1+d(-1)$ &
\multicolumn{1}{c|}{$J$} & \multicolumn{1}{c}{$a$} \\
\hline
O    & & & $1$ \\
I    & $1$ & $\{(7, 3)\}$ or $\{(3, 1), (5, 2)\}$ & $2$ \\
IIa  & $2$ & $\{(r, 2)\}$                         & $4/rE^3 = 2$ or $4$ \\
IIb  & $2$ & $\{(r_1, 1), (r_2, 1)\}$ $(r_1 \le r_2)$ &
$(r_1+r_2)/r_1r_2E^3 \ge2$ \\
III  & $3$ & $\{(r, 1)\}$                         & $(1+r)/rE^3 \ge2$ \\
IV   & $4$ & $\emptyset$                          & $2$ \\
\hline
\end{tabular}
\end{center}
\end{Theorem}

\begin{Remark}\label{rem:typeIV}
$f$ is of type IV if and only if
$f$ is the usual blowup along a smooth point $P$.
\end{Remark}

\begin{Remark}\label{rem:precise}
Type IIb is divided into types ${\textrm{IIb}}^{\vee}$
and ${\textrm{IIb}}^{\vee \vee}$ (\ref{rem:divide}).
By (\ref{rem:basket}) if $f$ is of type ${\textrm{IIb}}^{\vee}$ then
$J=\{(r,1), (r,1)\}$ or $\{(2,1), (4,1)\}$, and $a=2$, $3$ respectively.
In this case from \cite{Mo85} the unique non-Gorenstein point
$Q \in Y$ can be described as one of the following equations $\phi=0$
using semi-invariant local coordinates $y_1, y_2, y_3, y_4$
of the index one cover $Q^\sharp \in Y^\sharp$
with weights $w=\wt(y_1, y_2, y_3, y_4)$.

\begin{enumerate}
\item $J=\{(r, 1), (r, 1)\}$.
\begin{enumerate}
\item $\phi=y_1y_2+g(y_3^r, y_4)$ and $w=(1, -1, 2, 0)$.
The order of $g(0, y_4)$ is $2$.
The type of a general Du Val section of a germ $Q \in Y$ is
$A_{2r-1}$.\label{itm:A/r}
\item $r=3$, $\phi=y_4^2+\phi_3(y_1, y_2, y_3)+ \phi_{\ge 4}(y_1, y_2, y_3)$
and $w=(1, -1, 2, 0)$.
$\phi_3$ is $y_1^3+y_2^3+y_3^3$, $y_1^3+y_2y_3^2$ or $y_1^3+y_2^3$,
and $\phi_{\ge 4}$ is of order $\ge 4$.
The type of a general Du Val section of a germ $Q \in Y$ is
$E_6$.\label{itm:D/3}
\end{enumerate}
\item $J=\{(2, 1), (4, 1)\}$.
$\phi=y_1^2+y_2^2+g(y_3^2, y_4)$ and
$w=(1, 3, 3, 2)$.
The order of $g(0, y_4)$ is $3$.
The type of a general Du Val section of a germ $Q \in Y$ is
$D_5$.\label{itm:A/4}
\end{enumerate}
\end{Remark}

\begin{Remark}
By the proof of this theorem in \cite{Ka01},
we can restrict possible values of $J$ even if $f$ is of type $O$.
In fact $J$ is one of\\
$\{(7, 3)\}$, $\{(8, 3)\}$,\\
$\{(2, 1), (5, 2)\}$, $\{(3, 1), (5, 2)\}$, $\{(4, 1), (5, 2)\}$,
$\{(2, 1), (7, 2)\}$,\\
$\{(2, 1), (2, 1), (r_3, 1)\}$, $\{(2, 1), (3, 1), (3, 1)\}$,\\
$\{(2, 1), (3, 1), (4, 1)\}$, $\{(2, 1), (3, 1), (5, 1)\}$,\\
$\{(r, 2)\}$, $\{(r_1, 1), (r_2, 1)\}$, $\{(r, 1)\}$ and $\emptyset$.\\
In particular the number of fictitious singularities is at most three.
\end{Remark}

\begin{Example}\label{exl:3sing}
The weighted blowup of the c$D_4$ singularity
$o \in (x_1^2+x_2^3+x_3^3+x_4^6=0) \subset \bC^4$ with weights
$\wt(x_1, x_2, x_3, x_4) = (3, 2, 2, 1)$ in \cite[Example 4.6]{Ka01}
is an example whose $J$ is $\{(2, 1), (2, 1), (2, 1)\}$.
\end{Example}

\begin{Text}\label{txt:Q}
(\ref{thm:num.classn}) is obtained from the singular Riemann-Roch formula
\cite[Theorem 10.2]{Re87} and a relative vanishing theorem
\cite[Theorem 1-2-5]{KMM87} with the exact sequences
\begin{align}
0 \to \cO_Y((i-1)E) \to \cO_Y(iE) \to \cQ_i \to 0.
\end{align}
$\cQ_i$ are $\mathrm{S}_2$ by \cite[Proposition 5.26]{KoMo98}, and are
reflexive because they are locally free on the restriction $E^o \subseteq E$ of
the Gorenstein locus of $Y$. Thus, since $E$ is Cohen-Macaulay,
$\cQ_i$ are concretely given by
$i_*(\cO_Y(iE) \otimes \cO_E |_{E^o})$, where $i$ is the induced map
$E^o \hookrightarrow E$.

The singular Riemann-Roch implies
\begin{align}\label{eqn:A_i}
\chi(\cQ_i) &= \frac{1}{12} \{2(3i^2-3i+1)- 3(2i-1)a + a^2\}E^3\\
\nonumber &\quad + \frac{1}{12} E\cdot c_2(Y) + A_i - A_{i-1},\\
\nonumber A_i &=\sum_{Q\in I}
\Bigl( - \overline{ie} \frac{r_Q^2-1}{12r_Q} +
\sum_{j=1}^{\overline{ie}-1}
\frac{\overline{jb_Q}(r_Q-\overline{jb_Q})}{2r_Q} \Bigr).
\end{align}
On the other hand a relative vanishing theorem implies
\begin{align}\label{eqn:vanish}
R^jf_*\cO_Y(iE) = 0 \qquad (i \le a, j \ge 1).
\end{align}
\end{Text}

We summarise the formulae in \cite{Ka01} obtained from the above.

\begin{Proposition}\label{prp:formulae}
\begin{enumerate}
\item $rE^3 \in \bZ_{>0}$.\label{itm:rE^3}
\item $d(i)=\chi(\cQ_i)$ \textup{(}$i \le a$\textup{)}.
\item If $f$ is of type IIb or III and we set
$r_1=1$, $r_2=r$ when $f$ is of type III,
then
\begin{align*}
d(-i)=1+ \Bigl\lfloor \frac{i}{r_1} \Bigr\rfloor
\qquad (0\le i < \min\{r_2, a\}).
\end{align*}\label{itm:d(-i)}
\item Set $\displaystyle B_i=\sum_{Q\in I}
\frac{\overline{iv_Q}(r_Q-\overline{iv_Q})}{2r_Q}$. Then
\begin{align*}
\chi(\cQ_{-i})-\chi(\cQ_{i+1})=
\Bigl(i+\frac{1}{2}\Bigr)aE^3+B_{i+1}-B_i.
\end{align*}\label{itm:diffchi}
\end{enumerate}
\end{Proposition}

\begin{Remark}\label{rem:di}
The values of $d(-i)=\chi(\cQ_{-i})$ below are used later.
\begin{enumerate}
\item (type IIa) By direct calculation using (\ref{eqn:A_i}), we have
$d(-1)=1$, $d(-2)\ge2$ if $a=2$ or $(a, r)=(4, 5)$.\label{itm:diIIa}
\item (types IIb and III) $d(-i)=1+ \lfloor \frac{i}{r_1} \rfloor$
($1\le i < \min\{r_2, a\}$) (\ref{prp:formulae}(\ref{itm:d(-i)})),
and $d(-a)=2+ \lfloor \frac{a}{r_1} \rfloor$ if $a < r_2$ by
(\ref{lem:Q.a+1}) and the proofs of
\cite[Proposition 4.4.3]{Ka01} and
the latter part of \cite[Theorem 4.5]{Ka01}.\label{itm:diIIbIII}
\end{enumerate}
\end{Remark}

\begin{Lemma}\label{lem:Q.a+1}
$\cQ_{a+1} = \omega_E$ and $h^0(\omega_E)=h^1(\omega_E)=0$, $h^2(\omega_E)=1$.
\end{Lemma}

\begin{proof}
Let $E^o$ on $E$ be the restriction of the Gorenstein locus on $Y$
and let $i \colon E^o \hookrightarrow E$ be the induced map.
$\omega_E = i_* \omega_{E^o}$ since $\omega_E$ is $\mathrm{S}_2$.
$\cQ_{a+1} = \omega_E$ comes from $\cQ_{a+1} |_{E^o} = \omega_{E^o}$

$H^i(\omega_E)$ is the dual of $\Ext_E^{2-i}(\omega_E, \omega_E)$.
Therefore $h^0(\omega_E)=h^1(\omega_E)=0$ since $E$ is Cohen-Macaulay.
$h^2(\omega_E)=\dim\Hom(\omega_E, \omega_E)$.
Since $\omega_E$ is $\mathrm{S}_2$, the map
$\Hom(\omega_E, \omega_E) \to \Hom(\omega_{E^o}, \omega_{E^o})$ is injective.
On the other hand $\dim\Hom(\omega_{E^o}, \omega_{E^o})=h^0(\cO_{E^o})=1$.
Because $\Hom(\omega_E, \omega_E) \neq 0$, we obtain $h^2(\omega_E)=1$.
\end{proof}

\begin{Text}\label{txt:easy.elephant}
Let $H_X \ni P$ be a general hyperplane section on $X$
and let $H$ be its strict transform on $Y$.
Write $f^*H_X = H + bE$. $b$ is the largest integer satisfying that
$f_*\cO_Y(-bE) = \fm_P$.
Let $Y^o$ be the Gorenstein locus of $Y$. By the adjunction formula we have
\begin{align*}
\omega_H|_{Y^o \cap H} = f^*\omega_{H_X} \otimes \cO_Y((a-b)E)|_{Y^o \cap H}.
\end{align*}
$a \ge b$ because $H_X$ is canonical.
Moreover if $a=b$ then $H$ is normal and gives a Du Val section
in the anticanonical system of $Y$.
Combining it with (\ref{thm:num.classn}) we have the following.
\end{Text}

\begin{Theorem}\label{thm:O.I.elephant}
If $f$ is of type O or I, then the strict transform $H$ of a general hyperplane
section $H_X$ on $X$ gives a Du Val section in the anticanonical system on $Y$.
\end{Theorem}

\begin{Text}\label{txt:comm.diag}
In the remainder of this section we investigate the scheme $H \cap E$
of dimension $1$. We use the commutative diagrams below repeatedly.
\begin{align*}
\begin{array}{ccccccccc}
&&\cO_Y(iE) \otimes \cO_Y(-E)&\to&\cO_Y(iE)&\to&\cO_Y(iE) \otimes \cO_E&\to&0\\
&&\downarrow&&\|&&\downarrow&&\\
0&\to&\cO_Y((i-1)E)&\to&\cO_Y(iE)&\to&\cQ_i&\to&0,\\
&&&&&&\downarrow&&\\
&&&&&&0&&
\end{array}\\
\begin{array}{ccccccc}
0&&&&&\\
\uparrow&&&&&\\
\cO_Y(iE) \otimes \cI_{H \cap E \subset E}
&\!\!\to\!\!&\cO_Y(iE) \otimes \cO_E
&\!\!\to\!\!&\cO_Y(iE) \otimes \cO_{H \cap E}&\to0\\
\uparrow&&\|&&\|&\\
\cO_Y(iE) \otimes \cO_Y(-H) \otimes \cO_E
&\!\!\to\!\!&\cO_Y(iE) \otimes \cO_E
&\!\!\to\!\!&\cO_Y(iE) \otimes \cO_{H \cap E}&\to0\\
\downarrow&&\downarrow&&\downarrow&\\
0\to\quad\cQ_{i+b}\qquad\ &\!\!\to\!\!&\cQ_i&\!\!\to\!\!
&\cR_i&\to0,\\
&&\downarrow&&\downarrow&\\
&&0&&0&
\end{array}
\end{align*}
where $b$ is that in (\ref{txt:easy.elephant}).
\end{Text}

Chasing the diagrams (\ref{txt:comm.diag}) with (\ref{eqn:vanish}) and
(\ref{lem:Q.a+1}), we have the following.

\begin{Lemma}\label{lem:tree}
\begin{enumerate}
\item $h^j(\cO_Y(iE)\otimes \cO_E)=0$
      \textup{(}$i \le a, j \ge 1$\textup{)}.\\
      $h^1(\cO_Y((a+1)E)\otimes \cO_E)=0$,
      $h^2(\cO_Y((a+1)E)\otimes \cO_E)=1$.\label{itm:van}
\item If $f$ is not of type I, then\\
      $h^1(\cO_Y(iE)\otimes \cO_{H \cap E})=0$
      \textup{(}$i \le a-1$\textup{)}.\\
      $h^1(\cO_Y(aE)\otimes \cO_{H \cap E})=1$.\label{itm:curvevan}
\end{enumerate}
\end{Lemma}

\begin{Remark}\label{rem:tree}
$h^0(\cO_{H \cap E})=1$ and $h^1(\cO_{H \cap E})=0$ in any case,
even if $f$ is of type I (\ref{thm:O.I.elephant}).
Especially $H \cap E$ has no embedded points and
$(H \cap E)_{\red}$ is a union of $\bP^1$.
\end{Remark}

\begin{Text}\label{txt:normal.form}
The main ingredient we study is an irreducible reduced subscheme
$C \cong \bP^1 \subseteq H \cap E$.
Here we introduce a normal form of $Q \in C \subset Y$
at each non-Gorenstein point $Q$ of $Y$ through which $C$ passes,
following \cite[Lemma 2.7]{Mo88}.

Let $r_Q$ be the local index of $Q \in Y$.
Take the index one cover $\pi \colon (Y^\sharp \ni Q^\sharp) \to (Y \ni Q)$
and set $C^\sharp=(C \times_Y Y^\sharp)_{\red}$.
We have an identification (\ref{thm:non.Gor})
\begin{align*}
Q^\sharp \in Y^\sharp \cong o \in (\phi=0) \subset \bC^4=
x_1x_2x_3x_4\textrm{-}\mathrm{space}.
\end{align*}
Let $s_Q$ be the number of irreducible components of $C^\sharp$ and
$Q^\dag \in C^\dag$ be the normalisation
of one of the irreducible components of $C^\sharp$.
Let $t \in \cO_{C, Q}$ and $t^{s_Q/r_Q} \in \cO_{C^\dag, Q^\dag}$ be
uniformising parameters of $C$ and $C^\dag$.
Let $a_i$ be the minimal number such that there exists a semi-invariant
function with weight $\wt x_i$ whose image in $\cO_{C^\dag, Q^\dag}$ has
order $a_i/r_Q$ with respect to $t$.
We note that
\begin{align}\label{eqn:weights}
(a_1, a_2, a_3, a_4)
\in s_Q\bZ(\wt x_1, \wt x_2, \wt x_3, \wt x_4) \quad
\textrm{in } (\bZ/(r_Q))^4.
\end{align}
Then we can take an identification (\ref{thm:non.Gor})
such that $x_i |_{C^\dag}=t^{a_i/r_Q}$ ($i=1, 2, 3$), and
$x_4 |_{C^\dag}=t^{a_4/r_Q}$ if $Q^\sharp \in Y^\sharp$ is singular.
\end{Text}

\begin{Text}\label{txt:how.to.comp}
We give how to compute the images of some natural maps of sheaves
defined on a germ $Q \in C \subset Y$ using the data
$a_1, a_2, a_3, a_4$ in (\ref{txt:normal.form}).
For $n \in \bZ/(r_Q)$, we define $w_Q^C(n)$ as the smallest nonnegative integer
such that $(w_Q^C(n), n) \in \bZ \times \bZ/(r_Q)$ is contained in
the semigroup
\begin{align}\label{eqn:semigroup}
&\bZ_{\ge0}(a_1, \wt x_1)+\bZ_{\ge0}(a_2, \wt x_2)+
\bZ_{\ge0}(a_3, \wt x_3)+\bZ_{\ge0}(a_4, \wt x_4),\\
\nonumber\textrm{or } &\bZ_{\ge0}(a_1, \wt x_1)+\bZ_{\ge0}(a_2, \wt x_2)+
\bZ_{\ge0}(a_3, \wt x_3) \textrm{ if $Y^\sharp$ smooth}.
\end{align}
$w_Q^C(0)=0$, and $(r_Q, 0)$ is contained in the above semigroup
because $C$ is smooth.
Take a reflexive sheaf $\cL$ on $Y$ which is isomorphic to the ideal sheaf
defined by $x_1=0$ outside $Q$. We note that
$\cL^{[\otimes v_Q]} \cong \cO_Y(E)$ and
$\cL^{[\otimes b_Q]} \cong \cO_Y(K_Y)$ on the germ at $Q$,
where $\cL^{[\otimes i]}$ denotes the double dual of $\cL^{\otimes i}$.
Then for any integers $j_1, \ldots, j_k$ with $\sum_{1 \le i \le k} j_i=0$,
the image of the natural map
\begin{align*}
\cL^{[\otimes j_1]} \otimes \cdots \otimes
\cL^{[\otimes j_k]} \otimes \cO_C \to \cO_C
\end{align*}
is $(\fm_{Q \subset C})^{\sum_{1 \le i \le k} w_Q^C(-j_i)/r_Q}$, where
$\fm_{Q \subset C}$ is the ideal sheaf of $Q$ in $C$.
\end{Text}

\begin{Text}
We also need to treat sheaves $\cI_C/\cI_C^{(2)}$ as in \cite[Section 2]{Mo88}.
For any $C \cong \bP^1 \subset H \cap E$
its ideal sheaf as the subscheme in $Y$
is denoted by $\cI_C$, and its symbolic $2$-power is denoted by $\cI_C^{(2)}$.
Let $C^{(2)} \subset Y$ be the closed subscheme defined by $\cI_C^{(2)}$.
$C^{(2)} \subseteq 2H \cap 2E$.
$H^1(\cO_{2H \cap 2E})=0$ unless $f$ is of type O or I as in (\ref{lem:tree}),
and then $H^1(\cO_{C^{(2)}})=0$.
By the exact sequence
\begin{align*}
0 \to \cI_C/\cI_C^{(2)} \to \cO_{C^{(2)}} \to \cO_C \to 0,
\end{align*}
we have
\begin{sLemma}\label{lem:I/I2}
$H^1(\cI_C/\cI_C^{(2)})=0$ unless $f$ is of type O or I.
\end{sLemma}
Consider the natural map
\begin{align}\label{eqn:map}
\bigwedge^2 \cI_C/\cI_C^{(2)} \otimes \cO_C(K_C) \to [\cO_Y(K_Y)]_C,
\end{align}
where $[\cL]_C$ denotes the torsionfree part of $\cL \otimes \cO_C$.
The length of the cokernel of this map is
$2+\deg[\cO_Y(K_Y)]_C-\deg\bigwedge^2 \cI_C/\cI_C^{(2)}$.
By (\ref{lem:I/I2}) we obtain
\begin{sLemma}\label{lem:cokernel}
The length of the cokernel of the map \textup{(\ref{eqn:map})} is at most
$4+\deg[\cO_Y(K_Y)]_C$ unless $f$ is of type O or I.
\end{sLemma}
The next lemma is obtained by following \cite[Corollary 2.15]{Mo88}
almost faithfully.
We remark that $(2, 2) \in \bZ \times \bZ/(4)$ is contained
in the semigroup (\ref{eqn:semigroup})
when $Q$ is of exceptional type (\ref{itm:excep}) in (\ref{thm:non.Gor})
because $C$ is smooth at $Q$.
\begin{sLemma}[{\cite[Corollary 2.15]{Mo88}}]\label{lem:exclude}
Let $Q \in C \subset Y$ be a singular point of $Y$.
Then the map \textup{(\ref{eqn:map})} is not surjective at $Q$.
\end{sLemma}
\end{Text}

The next lemma demonstrates computation of the length of the cokernel of
the map (\ref{eqn:map}) in the simplest case.

\begin{Lemma}\label{lem:comp.cok}
Let $Q \in Y$
be a germ of a three dimensional terminal quotient singularity
of index $r$ obtained by the quotient of $\bC^3=x_1x_2x_3$-space
with weights $\wt(x_1, x_2, x_3)=(1, -1, b)$, and
let $Q \in C \subset Y$ be a smooth curve with uniformising parameter $t$
given by $(x_1, x_2, x_3) |_C=(t^{c/r}, t^{1-c/r}, 0)$
for some $0 < c < r$.
Then the length of the cokernel of the map \textup{(\ref{eqn:map})}
\begin{align*}
\bigwedge^2 \cI_C/\cI_C^{(2)} \otimes \cO_C(K_C) \to [\cO_Y(K_Y)]_C
\end{align*}
is $\min\{c, r-c\}$.
\end{Lemma}

\begin{proof}
First we give an explicit description of this map using local
coordinates. The map $\cO_Y(K_Y) \otimes \cO_C \to [\cO_Y(K_Y)]_C \cong \cO_C$
can be locally given, but not canonically, by
\begin{align}\label{eqn:comp.map1}
\cO_Y(K_Y) \otimes \cO_C &\to \cO_C \\
\nonumber f(x_1, x_2, x_3) dx_1 \wedge dx_2 \wedge dx_3 &\mapsto
\frac{f(t^{c/r}, t^{1-c/r}, 0)}{t^{w_Q^C(-b)/r}}.
\end{align}
$\cO_C(K_C)$ is generated by $d(x_1x_2)$.
Hence the length of the cokernel of our map equals that of the composition of
\begin{align}\label{eqn:comp.map2}
\cI_C/\cI_C^{(2)} \times \cI_C/\cI_C^{(2)} &\to \cO_Y(K_Y) \otimes \cO_C \\
\nonumber f_1 \times f_2 &\mapsto df_1 \wedge df_2 \wedge d(x_1x_2)
\end{align}
and (\ref{eqn:comp.map1}).

It is easy to see that $\cI_C/\cI_C^{(2)}$ is generated by functions of form
$g_{-b}(x_1, x_2) \cdot x_3$ and $g_c(x_1, x_2) \cdot (x_1^{r-c}-x_2^c)$
where $g_{-b}$ and $g_c$ are semi-invariant functions of weights
$-b$ and $c$ respectively. Therefore we have only to consider
images of $g_{-b} dx_3 \wedge
d(g_c\cdot (x_1^{r-c}-x_2^c)) \wedge d(x_1x_2)$ by the map
(\ref{eqn:comp.map1}).
Hence it is enough to consider
$g_c=x_1^c$ and $x_2^{r-c}$. In these cases
\begin{align*}
&g_{-b} dx_3 \wedge
d(g_c\cdot (x_1^{r-c}-x_2^c)) \wedge d(x_1x_2)\\
&=
\begin{cases}
g_{-b}\cdot rx_1^r dx_1 \wedge dx_2 \wedge dx_3
& (g_c=x_1^c)\\
g_{-b}\cdot rx_2^r dx_1 \wedge dx_2 \wedge dx_3
& (g_c=x_2^{r-c}).
\end{cases}
\end{align*}
Since we can choose $g_{-b}$ so that
$g_{-b}(t^{c/r}, t^{1-c/r})=t^{w_Q^C(-b)/r}$, our length is the minimum
of the orders of
$x_1^r$ and $x_2^r$ with $(x_1, x_2)=(t^{c/r}, t^{1-c/r})$
with respect to $t$, which is $\min\{c, r-c\}$.
\end{proof}

\begin{Text}\label{txt:comp.cok}
We need to compute the length in the more complicate case below.
Let $Q \in Y$ be a germ of a three dimensional terminal quotient singularity
of index $r$ obtained by the quotient of $\bC^3=x_1x_2x_3$-space
with weights $\wt(x_1, x_2, x_3)=(1, -1, b)$, and
let $Q \in C \subset Y$ be a smooth curve with uniformising parameter $t$
given by $(x_1, x_2, x_3) |_C=(t^{a_1/r}, t^{a_2/r}, t^{a_3/r})$
for some $a_1, a_2, a_3$.
We can take an invariant monomial $g(x_1, x_2, x_3)$ such that
$g(t^{a_1/r}, t^{a_2/r}, t^{a_3/r})=t$.

$\cI_C/\cI_C^{(2)}$ is generated by functions of form
$x_1^{s_1}x_2^{s_2}x_3^{s_3}-x_1^{t_1}x_2^{t_2}x_3^{t_3}$
with $(s_1,s_2,s_3)\neq(t_1,t_2,t_3)$ such that
$s_1-s_2+bs_3 \equiv t_1-t_2+bt_3 \equiv 0$ modulo $r$
and $a_1s_1+a_2s_2+a_3s_3=a_1t_1+a_2t_2+a_3t_3$.
Our length is the minimum of the orders of the images of
$d(x_1^{s_1}x_2^{s_2}x_3^{s_3}-x_1^{t_1}x_2^{t_2}x_3^{t_3}) \wedge
d(x_1^{\bar{s}_1}x_2^{\bar{s}_2}x_3^{\bar{s}_3}
-x_1^{\bar{t}_1}x_2^{\bar{t}_2}x_3^{\bar{t}_3}) \wedge dg$
by the map
\begin{align*}
\cO_Y(K_Y) \otimes \cO_C &\to \cO_C \\
\nonumber f(x_1, x_2, x_3) dx_1 \wedge dx_2 \wedge dx_3 &\mapsto
\frac{f(t^{a_1/r}, t^{a_2/r}, t^{a_3/r})}{t^{w_Q^C(-b)/r}},
\end{align*}
where $x_1^{s_1}x_2^{s_2}x_3^{s_3}-x_1^{t_1}x_2^{t_2}x_3^{t_3}$ and
$x_1^{\bar{s}_1}x_2^{\bar{s}_2}x_3^{\bar{s}_3}
-x_1^{\bar{t}_1}x_2^{\bar{t}_2}x_3^{\bar{t}_3}$ satisfy the conditions
mentioned above.
The order for fixed
$x_1^{s_1}x_2^{s_2}x_3^{s_3}-x_1^{t_1}x_2^{t_2}x_3^{t_3}$ and
$x_1^{\bar{s}_1}x_2^{\bar{s}_2}x_3^{\bar{s}_3}
-x_1^{\bar{t}_1}x_2^{\bar{t}_2}x_3^{\bar{t}_3}$ is $+\infty$ or
\begin{multline*}
\frac{1}{r}\cdot
\{(a_1s_1+a_2s_2+a_3s_3)+(a_1\bar{s}_1+a_2\bar{s}_2+a_3\bar{s}_3)+r\\
-(a_1+a_2+a_3+w_Q^C(-b))\}.
\end{multline*}
In particular if we can choose $g=x_1x_2$ then it is $+\infty$ or
\begin{align*}
\frac{1}{r}\cdot
\{(a_1s_1+a_2s_2+a_3s_3)+(a_1\bar{s}_1+a_2\bar{s}_2+a_3\bar{s}_3)
-(a_3+w_Q^C(-b))\}.
\end{align*}
\end{Text}

\section{Local analysis at a non-Gorenstein point}\label{sec:local}
\begin{Text}\label{txt:D}
In this section we analyse the structure of
$C \cong \bP^1 \subseteq H \cap E$ under the assumption that
$f$ is of type II or III focusing on the germ at a non-Gorenstein point.
We can see the theorems in this section automatically
when $P$ is a smooth point
because of the explicit description \cite{Ka01}, and hence
we here impose an extra assumption, that is, $P$ is a singular point.
We have $(H \cdot C) \le (H \cdot [H \cap E])=E^3$.
By (\ref{rem:a<r}), (\ref{thm:num.classn}),
(\ref{prp:formulae}(\ref{itm:rE^3})) and (\ref{rem:tree})
the configuration of the $1$-cycle $[H \cap E]=\sum [C_i]$ is as follows.
\begin{enumerate}
\item (type IIa) $a=2$ or $4$ and $Y$ has one non-Gorenstein point $Q$,
through which any $C_i$ passes.
\begin{enumerate}
\item $a=2$ and $H \cap E \cong \bP^1$ scheme theoretically.\label{itm:2ir}
\item $a=2$ and $[H \cap E]=[\bP^1]+[\bP^1]$ is
reducible or nonreduced.\label{itm:2r}
\item $a=4$ and $H \cap E \cong \bP^1$ scheme theoretically.\label{itm:4ir}
\end{enumerate}\label{itm:DIIa}
\item (type ${\textrm{IIb}}^{\vee}$)
$J=\{(r,1), (r,1)\}$ or $\{(2,1), (4,1)\}$ and $a=2$, $3$ respectively.
$Y$ has one non-Gorenstein point $Q$.
$H \cap E \cong \bP^1$ scheme theoretically
and passes through $Q$.\label{itm:DIIb^}
\item (type ${\textrm{IIb}}^{\vee \vee}$)
$a \le (r_1+r_2)/2$ and $Y$ has two non-Gorenstein points $Q_1$
and $Q_2$ of indices $r_1$ and $r_2$.
Exactly one of the irreducible components of $(H \cap E)_{\red}$
passes through both $Q_1$ and $Q_2$,
and the rest pass through $Q_2$ only.\label{itm:DIIb^2}
\item (type III) $a \le (1+r)/2$ and $Y$ has one non-Gorenstein point $Q$,
through which any $C_i$ passes.\label{itm:DIII}
\end{enumerate}
\end{Text}

\begin{Text}\label{txt:s_C}
Let $s_C(i)$ be the integer such that
$[\cO_Y(iE)]_C \cong \cO_{\bP^1}(s_C(i))$.
The following are easy to see by (\ref{rem:a<r}),
(\ref{lem:tree}) and (\ref{lem:quot}).
\begin{enumerate}
\item Assume that $C$ comes from type IIa. Then
\begin{enumerate}
\item $s_C(1)= \cdots= s_C(a-1)=-1$.
\item $s_C(-i) \ge -1$ ($i \ge 1$).
\item $s_C(a)=-2$ in (\ref{txt:D}(\ref{itm:2ir}), (\ref{itm:4ir})), and
$s_C(a)=-2$ or $-1$ in (\ref{txt:D}(\ref{itm:2r})).
\end{enumerate}\label{itm:sIIa}
\item Assume that $C$ comes from
type ${\textrm{IIb}}^{\vee}$. Then
\begin{enumerate}
\item $s_C(1)= \cdots= s_C(a-1)=-1$. $s_C(a)=-2$.
\item $s_C(-i) \ge -1$ ($i \ge 1$).
\end{enumerate}\label{itm:sIIb^}
\item Assume that $C$ comes from
type ${\textrm{IIb}}^{\vee \vee}$ with $Q_1, Q_2 \in C$. Then
\begin{enumerate}
\item $s_C(1)= \cdots= s_C(a-1)=-1$.
\item $s_C(-i)=-1$ ($1 \le i \le a-1$, $r_1 \nmid i$).
\item $s_C(-r_1i)= 0$ or $-1$ ($1 \le r_1i \le a-1$).
\item $(s_C(a), s_C(-a))=(-2, 0)$, $(-2, -1)$ or $(-1, -1)$.
\end{enumerate}\label{itm:sIIb}
\item Assume that $C$ comes from type III. Then
\begin{enumerate}
\item $s_C(1)= \cdots= s_C(a-1)=-1$.
\item $s_C(-1)= \cdots= s_C(-(a-1))=0$.
\item $(s_C(a), s_C(-a))=(-2, 1)$, $(-2, 0)$ or $(-1, 0)$.
\end{enumerate}\label{itm:sIIbIII}
\end{enumerate}

For instance we consider the case (\ref{itm:sIIbIII}).
$s_C(i)<0$ if $i\ge1$ by the existence of the canonical map
$\cO_Y(iE)^{\otimes r} \otimes \cO_C \to \cO_Y(irE) \otimes \cO_C$,
and thus $s_C(1)= \cdots= s_C(a-1)=-1$, $s_C(a)=-2$ or $-1$ by
(\ref{lem:tree}). On the other hand (\ref{lem:quot}) with $\cL=\cO_Y(E)$
implies that $s_C(-1)=0$, and hence $s_C(-i)\ge0$ if $i\ge1$
by the existence of the canonical map
$\cO_Y(-E)^{\otimes i} \otimes \cO_C \to \cO_Y(-iE) \otimes \cO_C$.
Hence $s_C(-1)= \cdots= s_C(-(a-1))=0$ and
$(s_C(a), s_C(-a))=(-2, 1)$, $(-2, 0)$ or $(-1, 0)$
because the map $\cO_Y(-iE) \otimes \cO_Y(iE) \otimes \cO_C \to \cO_C$
is not surjective when $1\le i \le a<r$ (\ref{rem:a<r}).
\end{Text}

\begin{Lemma}\label{lem:quot}
Let $Q \in Y$ be a germ of a three dimensional terminal quotient singularity
of index $r$ obtained by the quotient of $\bC^3=x_1x_2x_3$-space
with weights $\wt(x_1, x_2, x_3)=(1, -1, b)$.
Let $\cL$ be a reflexive sheaf on $Y$ which is isomorphic to the ideal sheaf
defined by $x_1=0$ outside $Q$.
Then the image of the natural map $\cL \otimes \cL^{[-1]} \to \cO_Y$
is the ideal sheaf of $Q$ in $Y$.
\end{Lemma}

\begin{proof}
It is enough to show that any invariant monomial $\neq 1$ of $x_1, x_2, x_3$
decomposes into two semi-invariant monomials with weights $1$ and $-1$,
but it is trivial since $b$ is coprime to $r$.
\end{proof}

The following proposition stores more information on $s_C$.

\begin{Proposition}\label{prp:store}
Assume that $f$ is not of type O or I, and
take a curve $C \cong \bP^1 \subseteq H \cap E$.
Let $i$ be a positive integer, and let $T$, $T'$ be 
surfaces on $Y$ defined by general elements in $\cO_Y(-iE)$.
\begin{enumerate}
\item Assume that $\chi(\cQ_i) \ge 2$.
Then $(T \cap E)_{\red}$ has an irreducible component which
properly intersects $(T' \cap E)_{\red}$.\label{itm:storeX2}
\item Assume that $\chi(\cQ_j)=1$ \textup{(}$1 \le j < i$\textup{)}
and $\chi(\cQ_i) \ge 2$.
If $C=H \cap E$ scheme theoretically,
then $(T \cap E)_{\red}$ properly intersects $(T' \cap E)_{\red}$ and $C$.
Moreover $s_C(-j)=-1$ \textup{(}$1 \le j < i$\textup{)}
and $s_C(-i) \ge 0$.\label{itm:storeX1X2}
\end{enumerate}
\end{Proposition}

\begin{proof}
(\ref{itm:storeX2})
By (\ref{eqn:vanish}) we have a surjective map
$f_*\cO_Y(-iE) \twoheadrightarrow H^0(\cQ_{-i})$.
Write $E^o$ as the restriction of the Gorenstein locus of $Y$ to $E$.
Take a resolution $\hat{E}$ of $E$ and let $\hat{E^o}$ be the preimage
of $E^o$. We note that $H^0(\cQ_{-i} |_{E^o})=H^0(\cQ_{-i})$ (\ref{txt:Q}).
Consider a linear system $L$ on $\hat{E}$, not necessarily complete,
whose restriction to $\hat{E^o}$ is $H^0(\cQ_{-i} |_{E^o})$.
$h^0(\cQ_{-i})=\chi(\cQ_{-i}) \ge 2$ means that $L$ moves,
and its general member has an irreducible component which is not contained in
the union of the support of the fixed locus of $L$ and the preimage of
the nonnormal locus of $E$.
The element in $H^0(\cQ_{-i})$ corresponding to this member gives
a desired section $T$.

(\ref{itm:storeX1X2})
By (\ref{txt:comm.diag}) and (\ref{lem:tree}(\ref{itm:van}))
we have surjective maps
\begin{align*}
f_*\cO_Y(-jE)  &\twoheadrightarrow H^0(\cQ_{-j})
\twoheadrightarrow H^0(\cR_{-j})\twoheadrightarrow
H^0([\cO_Y(-jE)]_C) \quad (1 \le j \le i).
\end{align*}
We note that $E$ is smooth at the generic point of $C$
because $C=H \cap E$ scheme theoretically.
$H^0(\cQ_{-j})=\chi(\cQ_{-j}) =1$ ($1 \le j < i$) and
a nonzero element in it corresponds to, at the generic point of $C$,
the Cartier divisor defined by $jH$.
Hence the map $H^0(\cQ_{-j}) \twoheadrightarrow H^0([\cO_Y(-jE)]_C)$
is a zero map, and we obtain $s_C(-j)=-1$
by the above map and (\ref{lem:tree}(\ref{itm:curvevan})).

On the other hand the $L$ constructed from $H^0(\cQ_{-i})$
as in the proof of (\ref{itm:storeX2}) moves.
Because $h^0(\cQ_{-j})=1$ ($1 \le j < i$)
and $E$ is smooth at the generic point of $C$,
its general member does not have the pull-back of $C$ as its component.
Hence we have a desired section and we see that
the map $H^0(\cQ_{-i}) \twoheadrightarrow H^0([\cO_Y(-iE)]_C)$
is a nonzero map, whence $s_C(-i) \ge 0$.
\end{proof}

First we treat type IIa.

\begin{Theorem}\label{thm:localIIa}
Assume that $f$ is of type IIa.
Let $Q \in C \cong \bP^1 \subseteq H \cap E$ be the non-Gorenstein point
with an irreducible component of $H \cap E$
\textup{(\ref{txt:D}(\ref{itm:DIIa}))}.
Take a normal form of $Q \in C \subset Y$ \textup{(\ref{txt:normal.form})}.
\begin{enumerate}
\item The case \textup{(\ref{itm:2ir})} in \textup{(\ref{txt:D})}
does not happen.\label{itm:th2ir}
\item In the case \textup{(\ref{itm:2r})} in \textup{(\ref{txt:D})},
we have $s_C(-2)=0$ and $(a_1, a_2, a_3)=((r+1)/2, (r-1)/2 , 2)$.
In particular we may choose semi-invariant local coordinates
$x_1, x_2, x_3$ with weights $\wt(x_1, x_2, x_3)=(1, -1, 4)$
of the index one cover $Q^\sharp \in  Y^\sharp$ so that
$(x_1, x_2, x_3) |_{C^\dag}=(t^{(r+1)/2r}, 0, t^{2/r})$
\textup{(}$r \equiv 1$ modulo $4$\textup{)} or
$(0, t^{(r-1)/2r}, t^{2/r})$
\textup{(}$r \equiv 3$ modulo $4$\textup{)}.\label{itm:th2r}
\item In the case \textup{(\ref{itm:4ir})} in \textup{(\ref{txt:D})},
we have $r=5$, $s_C(-4)=0$ and $(a_1, a_2, a_3)=(3, 2, 4)$.
In particular we may choose semi-invariant local coordinates
$x_1, x_2, x_3$ with weights $\wt(x_1, x_2, x_3)=(1, -1, 3)$
of the index one cover $Q^\sharp \in  Y^\sharp$ so that
$(x_1, x_2, x_3) |_{C^\dag}=(t^{3/5},t^{2/5}, 0)$.\label{itm:th4ir}
\end{enumerate}
\end{Theorem}

\begin{proof}
(\ref{itm:th2ir}).
Compute the image of the map
\begin{align*}
\cO_Y(E)^{\otimes r} \otimes \cO_Y(-rE) \otimes \cO_C \to \cO_C
\end{align*}
by (\ref{txt:how.to.comp}). It is $(\fm_{Q \subset C})^{w_Q^C(-2)}$.
On the other hand this image is the same as that of
\begin{align*}
[\cO_Y(E)]_C^{\otimes r} \otimes [\cO_Y(-rE)]_C \to \cO_C,
\end{align*}
whence $rs_C(1)-r(E \cdot C)=-w_Q^C(-2)$.
By $s_C(1)=-1$ (\ref{txt:s_C}(\ref{itm:sIIa})) and $(H \cdot C)=E^3=2/r$
(\ref{thm:num.classn}) we have $w_Q^C(-2)=r-2$. Hence we can write
$(a_1, a_2, a_3)=(1+rm_1, r-1+rm_2, 4+rm_3)$ with
$m_1, m_2, m_3 \in \bZ_{\ge 0}$ by (\ref{eqn:weights}).
We see $m_1=0$ easily since there exists an invariant monomial of
$x_1, x_2, x_3$ whose restriction on $C$ is $t$.
Thus $(a_1, a_2, a_3)=(1, r-1, 4)$ and especially $w_Q^C(-4)=r-4$.

Consider another map
$\cO_Y(2E)^{\otimes r} \otimes \cO_Y(-2rE) \otimes \cO_C \to \cO_C$.
Then we have $rs_C(2)-2r(E \cdot C)=-w_Q^C(-4)$, whence
$s_C(2)=-1$ by $w_Q^C(-4)=r-4$ and $(H \cdot C)=2/r$.
But this contradicts (\ref{txt:s_C}(\ref{itm:sIIa})).

(\ref{itm:th2r})
By the map
$\cO_Y(E)^{\otimes r} \otimes \cO_Y(-rE) \otimes \cO_C \to \cO_C$,
we have $rs_C(1)-r(E \cdot C)=-w_Q^C(-2)$.
Hence $w_Q^C(-2)=r-1$ by (\ref{txt:s_C}(\ref{itm:sIIa})) and
$(H \cdot C)=E^3/2=1/r$ (\ref{thm:num.classn}),
and we can write $(a_1, a_2, a_3)=((r+1)/2+rm_1, (r-1)/2+rm_2, 2+rm_3)$
with $m_1, m_2, m_3 \in \bZ_{\ge 0}$ by (\ref{eqn:weights}).
Because there exists an invariant monomial of
$x_1, x_2, x_3$ whose restriction on $C$ is $t$,
we have one of the following.
\begin{enumerate}
\item $(r, 0) \in \bZ_{>0}(a_1, 1)+\bZ_{>0}(a_3, 4)$
or $\in \bZ_{>0}(a_2, -1)+\bZ_{>0}(a_3, 4)$
in $\bZ \times \bZ/(r)$.\label{itm:pp2}
\item $r=a_1+a_2$.\label{itm:pp3}
\end{enumerate}
If (\ref{itm:pp2}) holds then $(a_1, a_2, a_3)=((r+1)/2, (r-1)/2, 2)$
and $w_Q^C(4)=4$. By the map
$\cO_Y(-2E)^{\otimes r} \otimes \cO_Y(2rE) \otimes \cO_C \to \cO_C$,
we have $rs_C(-2)+2r(E \cdot C)=-w_Q^C(4)$, whence $s_C(-2)=0$.
Now we assume (\ref{itm:pp3}).
Then $(a_1, a_2, a_3)=((r+1)/2, (r-1)/2, 2+rm_3)$ with some $m_3$.
We have nothing to do if $m_3=0$. Thus we also assume that $m_3 > 0$
and $r \ge 7$. We note that $m_3 \le 2$ by $4a_1=2+2r$.

If $m_3=2$ then we may choose semi-invariant local coordinates
$x_1, x_2, x_3$ with weights $\wt(x_1, x_2, x_3)=(1, -1, 4)$
of $Q^\sharp \in  Y^\sharp$ so that
$(x_1, x_2, x_3) |_{C^\dag}=(t^{(r+1)/2r}, t^{(r-1)/2r}, 0)$.
We see $w_Q^C(-4)=4w_Q^C(-1)=2r-2$ from $r \ge 7$.
Considering the map
$\cO_Y(2E)^{\otimes r} \otimes \cO_Y(-2rE) \otimes \cO_C \to \cO_C$,
we obtain $rs_C(2)-2r(E \cdot C)=-w_Q^C(-4)$, whence $s_C(2)=-2$.
By (\ref{lem:cokernel}), the length of the cokernel of the
the map (\ref{eqn:map}) at $Q$ is at most $2$.
Hence we have $(r-1)/2 \le 2$ by (\ref{lem:comp.cok}),
which contradicts $r \ge 7$.

Now it remains to exclude
the case where $(a_1, a_2, a_3)=((r+1)/2, (r-1)/2, r+2)$ with $r \ge 7$.
But if $r=7$ then $(a_1, a_2, a_3)=(4, 3, 9)$ and
we may choose the coordinates
so that $(x_1, x_2, x_3) |_{C^\dag}=(t^{4/7}, t^{3/7}, 0)$.
Then $s_C(2)=-2$ follows and contradicts
(\ref{lem:cokernel}) and (\ref{lem:comp.cok}).
Therefore we may assume that $r \ge 9$.
By the map
$\cO_Y(-3E)^{\otimes r} \otimes \cO_Y(3rE) \otimes \cO_C \to \cO_C$,
we have $rs_C(-3)+3r(E \cdot C)=-w_Q^C(6)$, whence
$w_Q^C(6)=3$ or $r+3$ since $s_C(-3) \ge -1$.
From $(a_1, a_2, a_3)=((r+1)/2, (r-1)/2, r+2)$ with $r \ge 9$,
$r$ has to be $9$.

$(a_1, a_2, a_3)=(5, 4, 11)$ and $w_Q^C(-4)=16$ when $r=9$.
By the map
$\cO_Y(2E)^{\otimes r} \otimes \cO_Y(-18E) \otimes \cO_C \to \cO_C$,
we obtain $9s_C(2)-18(E \cdot C)=-w_Q^C(-4)$, whence $s_C(2)=-2$.
We consider $x_1^{s_1}x_2^{s_2}x_3^{s_3}-x_1^{t_1}x_2^{t_2}x_3^{t_3}$
with $(s_1,s_2,s_3)\neq(t_1,t_2,t_3)$ such that
$s_1-s_2+4s_3 \equiv t_1-t_2+4t_3 \equiv 0$ modulo $9$
and $5s_1+4s_2+11s_3=5t_1+4t_2+11t_3$.
We can easily see that $(5s_1+4s_2+11s_3)/9 \ge 3$ for any such
$x_1^{s_1}x_2^{s_2}x_3^{s_3}-x_1^{t_1}x_2^{t_2}x_3^{t_3}$.
Hence by (\ref{txt:comp.cok}), the length of the cokernel of the map
(\ref{eqn:map}) is greater than or equal to $3+3-(a_3+w_Q^C(-4))/9=3$,
which contradicts $s_C(2)=-2$ and (\ref{lem:cokernel}).

(\ref{itm:th4ir})
By the map
$\cO_Y(iE)^{\otimes r} \otimes \cO_Y(-irE) \otimes \cO_C \to \cO_C$
($i=1, 2, 3, 4$),
we have $rs_C(i)-ir(E \cdot C)=-w_Q^C(-2i)$, whence
$w_Q^C(-2i)=r-i$ ($i=1, 2, 3$) and $w_Q^C(-8)=2r-4$
(\ref{txt:s_C}(\ref{itm:4ir})).
We can write
$(a_1, a_2, a_3)=((r+1)/2+rm_1, (r-1)/2+rm_2, 4+rm_3)$ with
$m_1, m_2, m_3 \in \bZ_{\ge 0}$ by (\ref{eqn:weights}).
If $r=5$ then $a_1=w_Q^C(-4)=3$ and $a_2=w_Q^C(-6)=2$ and
we have the desired description.
In this case $w_Q^C(8)=4$, and by the map
$\cO_Y(-4E)^{\otimes r} \otimes \cO_Y(4rE) \otimes \cO_C \to \cO_C$,
we have $rs_C(-4)+4r(E \cdot C)=-w_Q^C(8)$, whence $s_C(-4)=0$.

Now we derive a contradiction supposing that $r \ge 7$.
Since $w_Q^C(-4)=r-2$, we can decompose $r-2$ into
$(r+1)/2+4c$ or $(r-1)/2+4c$ for some positive integers $c$,
and then $r \equiv 5$, $3$ modulo $8$ respectively.
Since $w_Q^C(-6)=r-3$, similarly we can decompose $r-3$ into
$(r+1)/2+4c$, $(r-1)/2+4c$ or $4c$
and then $r \equiv 7$, $5$ modulo $8$, $r \equiv 3$ modulo $4$ respectively.
Hence $r \equiv 3$ modulo $8$ with $m_2=m_3=0$, or
$r \equiv 5$ modulo $8$ with $m_1=m_2=m_3=0$.
In particular $r \ge 11$ and $a_3=4$.

We consider $x_1^{s_1}x_2^{s_2}x_3^{s_3}-x_1^{t_1}x_2^{t_2}x_3^{t_3}$
with $(s_1,s_2,s_3)\neq(t_1,t_2,t_3)$ such that
$s_1-s_2+8s_3 \equiv t_1-t_2+8t_3 \equiv 0$ modulo $r$
and $a_1s_1+a_2s_2+a_3s_3=a_1t_1+a_2t_2+a_3t_3$.
We can easily see that $(a_1s_1+a_2s_2+a_3s_3)/r=2$ only if
each of $x_1^{s_1}x_2^{s_2}x_3^{s_3}$ and $x_1^{t_1}x_2^{t_2}x_3^{t_3}$ is
$x_1^2x_2^2$, $x_1x_3^{(3r-1)/8}$ ($r \equiv 3$ modulo $8$) or
$x_2x_3^{(3r+1)/8}$ ($r \equiv 5$ modulo $8$) and
otherwise $(a_1s_1+a_2s_2+a_3s_3)/r \ge 3$.
Hence by (\ref{txt:comp.cok}), the length of the cokernel of the map
(\ref{eqn:map}) is greater than or equal to $2+3-(a_3+w_Q^C(-8))/r=3$,
which contradicts $s_C(4)=-2$ (\ref{txt:s_C}(\ref{itm:sIIa}))
and (\ref{lem:cokernel}).
\end{proof}

It is easy to see the local structure in type ${\textrm{IIb}}^{\vee}$.

\begin{Theorem}\label{thm:localIIb^}
Assume that $f$ is of type ${\textrm{IIb}}^{\vee}$.
Let $Q \in C = H \cap E \cong \bP^1$ be the non-Gorenstein point
\textup{(\ref{txt:D}(\ref{itm:DIIb^}))}.
Take a normal form of $Q \in C \subset Y$ \textup{(\ref{txt:normal.form})}.
\begin{enumerate}
\item If $J=\{(r, 1), (r, 1)\}$, then $a=2$, $s_C(-2)=0$ and
$(a_1, a_2, a_3, a_4)=(r+1, r-1, 2, r)$.
In particular we may choose semi-invariant local coordinates
$x_1, x_2, x_3, x_4$ with weights $\wt(x_1, x_2, x_3, x_4)=(1, -1, 2, 0)$
of the index one cover $Q^\sharp \in  Y^\sharp$ so that
$(x_1, x_2, x_3, x_4) |_{C^\dag}=(0, 0, t^{2/r}, t)$.\label{itm:dIIb^r}
\item If $J=\{(2, 1), (4, 1)\}$, then $a=3$, $s_C(-3)=0$ and
$(a_1, a_2, a_3, a_4)=(5, 3, 3, 2)$.
In particular we may choose semi-invariant local coordinates
$x_1, x_2, x_3, x_4$ with weights $\wt(x_1, x_2, x_3, x_4)=(1, 3, 3, 2)$
of the index one cover $Q^\sharp \in  Y^\sharp$ so that
$(x_1, x_2, x_3, x_4) |_{C^\dag}=(0, 0, t^{3/4}, t^{1/2})$.\label{itm:dIIb^4}
\end{enumerate}
\end{Theorem}

\begin{proof}
(\ref{itm:dIIb^r})
First applying (\ref{prp:store}(\ref{itm:storeX1X2})) with $i=2$
(\ref{rem:di}(\ref{itm:diIIbIII})) we have $s_C(-1)=-1$ and $s_C(-2)\ge 0$.
The map
$\cO_Y(-2E) \otimes \cO_Y(E) \otimes \cO_C \to \cO_Y(-E) \otimes \cO_C$
implies that $s_C(-2)+s_C(1) \le s_C(-1)$. Hence $s_C(-2)=0$.
By the map
$\cO_Y(-E)^{\otimes r} \otimes \cO_Y(rE) \otimes \cO_C \to \cO_C$,
we have $rs_C(-1)+r(E \cdot C)=-w_Q^C(1)$.
By $s_C(-1)=-1$ and $E^3=1/r$ we have $w_Q^C(1)=1+r$. Hence we can write
$(a_1, a_2, a_3, a_4)=(r+1+rm_1, r-1+rm_2, 2+rm_3, r)$ with
$m_1, m_2, m_3 \in \bZ_{\ge 0}$ by (\ref{eqn:weights}) and the property
that $(r, 0)$ is contained in the semigroup (\ref{eqn:semigroup}).
Consider another map
$\cO_Y(-2E)^{\otimes r} \otimes \cO_Y(2rE) \otimes \cO_C \to \cO_C$.
Then we have $rs_C(-2)+2r(E \cdot C)=-w_Q^C(2)$.
Therefore $w_Q^C(2)=2$ from $s_C(-2)=0$,
whence $(a_1, a_2, a_3, a_4)=(r+1, r-1, 2, r)$.

(\ref{itm:dIIb^4})
By (\ref{prp:store}(\ref{itm:storeX1X2})) with large $i$
(\ref{rem:di}(\ref{itm:diIIbIII}))) we have $s_C(-1)=-1$.
By the map
$\cO_Y(-E)^{\otimes 4} \otimes \cO_Y(4E) \otimes \cO_C \to \cO_C$,
we have $4s_C(-1)+4(E \cdot C)=-w_Q^C(1)$.
By $s_C(-1)=-1$ and $E^3=1/4$ we obtain $w_Q^C(1)=5$. Hence we can write
$(a_1, a_2, a_3, a_4)=(5+4m_1, 3+4m_2, 3+4m_3, 2)$ with
$m_1, m_2, m_3 \in \bZ_{\ge 0}$.
Consider another map
$\cO_Y(E)^{\otimes 4} \otimes \cO_Y(-4E) \otimes \cO_C \to \cO_C$.
Then we have $4s_C(1)-4(E \cdot C)=-w_Q^C(-1)$.
Thus $w_Q^C(-1)=3$ (\ref{txt:s_C}(\ref{itm:sIIb^})),
whence $(a_1, a_2, a_3, a_4)=(5, 3, 3, 2)$.
$w_Q^C(3)=3$. By the map
$\cO_Y(-3E)^{\otimes 4} \otimes \cO_Y(12E) \otimes \cO_C \to \cO_C$,
we have $4s_C(-3)+12(E \cdot C)=-w_Q^C(3)$, whence $s_C(-3)=0$.
\end{proof}

\begin{Text}\label{txt:normIIb^2III}
It is hardest to treat types ${\textrm{IIb}}^{\vee \vee}$ and III.
\begin{enumerate}
\item If $f$ is of type ${\textrm{IIb}}^{\vee \vee}$,
then there exists the unique irreducible component
$C \cong \bP^1 \subseteq H \cap E$ through $Q_1$ and $Q_2$.
We do not consider any other components from now on.
Take normal forms of $Q_i \in C \subset Y$
\textup{(}$i=1, 2$\textup{)} using local coordinates $x_{i1}, x_{i2}, x_{i3}$
and numbers $a_{i1}, a_{i2}, a_{i3}$ respectively
(\ref{txt:normal.form}).\label{itm:normIIb^2}
\item If $f$ is of type III,
then we take any irreducible component $C \cong \bP^1 \subseteq H \cap E$
and a normal form of $Q \in C \subset Y$
(\ref{txt:normal.form}).\label{itm:normIII}
\end{enumerate}
\end{Text}

\begin{Proposition}\label{prp:tfae}
Assume that
$s_C(-a)=0$ in \textup{(\ref{txt:normIIb^2III}(\ref{itm:normIIb^2}))}
or $s_C(-a)=1$ in \textup{(\ref{txt:normIIb^2III}(\ref{itm:normIII}))}.
Then $C=H \cap E$ scheme theoretically, $s_C(a)=-2$, and
for any non-Gorenstein point $Q \in Y$ with
$Q \in C \cong \bP^1 =H \cap E$ in \textup{(\ref{txt:normIIb^2III})}
we may choose semi-invariant local coordinates
$x_1, x_2, x_3$ with weights
$\wt(x_1, x_2, x_3)=(1, -1, a)$
of the index one cover $Q^\sharp \in Y^\sharp$ so that
\begin{align*}
Q^\sharp \in C^\sharp \subset Y^\sharp
\cong o \in (x_3\textrm{-}\mathrm{axis})
\subset x_1x_2x_3\textrm{-}\mathrm{space}.
\end{align*}
\end{Proposition}

\begin{proof}
By the map
$\cO_Y(-aE)^{\otimes r} \otimes \cO_Y(arE) \otimes \cO_C \to \cO_C$,
we have
\begin{align*}
rs_C(-a)+ar(E \cdot C)=
\begin{cases}
-rw_{Q_1}^C(a)/r_1-rw_{Q_2}^C(a)/r_2
& (\ref{txt:normIIb^2III}(\ref{itm:normIIb^2}))\\
-w_Q^C(a) & (\ref{txt:normIIb^2III}(\ref{itm:normIII})).
\end{cases}
\end{align*}
Therefore by (\ref{thm:num.classn}),
\begin{align*}
\frac{w_{Q_1}^C(a)}{r_1}+\frac{w_{Q_2}^C(a)}{r_2} &= a(H \cdot C) \le
\frac{1}{r_1}+\frac{1}{r_2} &(\ref{txt:normIIb^2III}(\ref{itm:normIIb^2})),\\
\frac{w_Q^C(a)}{r} &= a(H \cdot C)-1 \le
\frac{1}{r} & (\ref{txt:normIIb^2III}(\ref{itm:normIII})).
\end{align*}
Hence $C=H \cap E$ (\ref{rem:tree}) and $w_Q^C(a)=1$ for any $Q \in C$
in (\ref{txt:normIIb^2III}), and we have the desired description.
$s_C(a)=-2$ follows from (\ref{lem:tree}(\ref{itm:curvevan})).
\end{proof}

\begin{Proposition}\label{prp:gen.by.two}
Assume that $f$ is of type ${\textrm{IIb}}^{\vee \vee}$ or III.
Consider a non-Gorenstein point $Q \in Y$ with
$Q \in C \cong \bP^1 \subseteq H \cap E$
given in \textup{(\ref{txt:normIIb^2III})}.
Then one of the following holds.
\begin{enumerate}
\item One of $(a_1, 1)$, $(a_2, -1)$
is generated by the other and $(a_3, a)$
in $\bZ \times \bZ/(r_Q)$,
and we may choose semi-invariant local coordinates
$x_1, x_2, x_3$ with weight $\wt(x_1, x_2, x_3)=(1, -1, a)$
of the index one cover $Q^\sharp \in Y^\sharp$ so that
$(x_1, x_2, x_3) |_{C^\dag}=(0, t^{c_Q/r_Q}, t^{1-\overline{ac_Q}/r_Q})$
or $(t^{c_Q/r_Q}, 0, t^{\overline{ac_Q}/r_Q})$
for some $0<c_Q<r_Q$.\label{itm:gen.two.ord}
\item $f$ is of type III,
$a=2$, $r \ge 5$, $s_Q(2)=-2$ and $(a_1, a_2, a_3)=(2, r-2, 4)$
or $(r-2, 2, 2r-4)$.
In particular we may choose semi-invariant local coordinates
$x_1, x_2, x_3$ with weights $\wt(x_1, x_2, x_3)=(1, -1, 2)$
of the index one cover $Q^\sharp \in  Y^\sharp$ so that
$(x_1, x_2, x_3) |_{C^\dag}=(t^{2/r}, t^{1-2/r}, 0)$ or
$(t^{1-2/r}, t^{2/r}, 0)$.
The value of $(H \cdot E)$ is $2/r$, $1-2/r$
respectively.\label{itm:gen.two.exc}
\end{enumerate}
\end{Proposition}

\begin{proof}
Recall that there exists an invariant monomial of $x_1, x_2, x_3$
whose restriction on $C$ is $t$.
Thus one of the following holds.
\begin{enumerate}
\item one of $a_1, a_2, a_3$ is $1$.\label{itm:p1}
\item $(r_Q, 0) \in \bZ_{>0}(a_1, 1)+\bZ_{>0}(a_3, a)$
or $\in \bZ_{>0}(a_2, -1)+\bZ_{>0}(a_3, a)$
in $\bZ \times \bZ/(r_Q)$.\label{itm:p2}
\item $r_Q=a_1+a_2$.\label{itm:p3}
\end{enumerate}
If (\ref{itm:p1}) or (\ref{itm:p2}) happens, then
one of $(a_1, 1)$ and $(a_2, -1)$
is generated by the other and $(a_3, a)$.
Therefore we may assume that only (\ref{itm:p3}) holds from now on.
By permutation of $x_i$ we may moreover assume that $r_Q \ge 5$.
Then it is easy to see that $w_Q^C(i)+w_Q^C(-i) \neq r_Q$ ($i=2, 3$),
whence the image of $[\cO_Y(-iE)]_C \otimes [\cO_Y(iE)]_C \to \cO_C$ at $Q$
is properly contained in $\fm_{Q \subset C}$ ($i=2, 3$)
(\ref{txt:how.to.comp}).
By (\ref{txt:s_C}) one of the following occurs.
\begin{enumerate}
\item $a=2$, $s_C(2)=-2$, $s_C(-2)=-1$ in
(\ref{txt:normIIb^2III}(\ref{itm:normIIb^2}))
or $s_C(-2)=0$ in (\ref{txt:normIIb^2III}(\ref{itm:normIII})).\label{itm:exc1}
\item $Q=Q_2$, $a=3$, $r_1=2$, $s_C(3)=-2$ and $s_C(-2)=s_C(-3)=-1$
in (\ref{txt:normIIb^2III}(\ref{itm:normIIb^2})).\label{itm:exc2}
\end{enumerate}

Consider the case (\ref{itm:exc1}).
$w_Q^C(1)+w_Q^C(-1)=r_Q$ and $w_Q^C(2)+w_Q^C(-2)=2r_Q$,
whence $w_Q^C(2)=2w_Q^C(1)$ and $w_Q^C(-2)=2w_Q^C(-1)$.
We have $(a_1, a_2, a_3)=(c, r_Q-c, 2c)$ by putting $w_Q^C(1)=c$,
and thus we may choose semi-invariant local coordinates
$x_1, x_2, x_3$ with weights $\wt(x_1, x_2, x_3)=(1, -1, 2)$
of the index one cover $Q^\sharp \in  Y^\sharp$ so that
$(x_1, x_2, x_3) |_{C^\dag}=(t^{c/r_Q}, t^{1-c/r_Q}, 0)$.

By $s_C(2)=-2$, (\ref{lem:cokernel}) and (\ref{lem:exclude}),
the length of the map (\ref{eqn:map}) at $Q$ is at most $1$
in (\ref{txt:normIIb^2III}(\ref{itm:normIIb^2}))
and at most $2$ in (\ref{txt:normIIb^2III}(\ref{itm:normIIb^2})).
Hence we have the result by (\ref{lem:comp.cok})
except the value of $(H \cdot C)$.
This value can be calculated from $rs_C(-1)+r(E \cdot C)=-w_Q^C(1)$
obtained by the map
$\cO_Y(-E)^{\otimes r} \otimes \cO_Y(rE) \otimes \cO_C \to \cO_C$.

Consider the case (\ref{itm:exc2}). $w_Q^C(3)+w_Q^C(-3)=2r_Q$.
By $s_C(3)=-2$ and (\ref{lem:cokernel}) and (\ref{lem:exclude}),
the length of the map (\ref{eqn:map}) at $Q$ must be $1$.
On the other hand, consider
$x_1^{s_1}x_2^{s_2}x_3^{s_3}-x_1^{t_1}x_2^{t_2}x_3^{t_3}$
with $(s_1,s_2,s_3)\neq(t_1,t_2,t_3)$ such that
$s_1-s_2+3s_3 \equiv t_1-t_2+3t_3 \equiv 0$ modulo $r_Q$
and $a_1s_1+a_2s_2+a_3s_3=a_1t_1+a_2t_2+a_3t_3$.
We can easily see that $(a_1s_1+a_2s_2+a_3s_3)/r_Q\ge2$ for any such
$x_1^{s_1}x_2^{s_2}x_3^{s_3}-x_1^{t_1}x_2^{t_2}x_3^{t_3}$
because $x_1x_2$ is the only invariant monomial
whose restriction to $C^\dag$ is $t$.
Hence by (\ref{txt:comp.cok}), the length of the cokernel of the map
(\ref{eqn:map}) at $Q$ is greater than or equal to $2+2-(a_3+w_Q^C(-3))/r_Q=2$,
which is a contradiction.
\end{proof}

We need more accurate description developing (\ref{prp:gen.by.two}).

\begin{Lemma}\label{lem:one.one}
Consider the case in \textup{(\ref{prp:gen.by.two}(\ref{itm:gen.two.ord}))}.
\begin{enumerate}
\item If $f$ is of type ${\textrm{IIb}}^{\vee \vee}$
\textup{(\ref{txt:normIIb^2III}(\ref{itm:normIIb^2}))},
then we may choose the coordinates so that
the restriction to $C^\dag$ at one of $Q_1$ and $Q_2$ is
of form $(0, t^{c_i/r_i}, t^{1-\overline{ac_i}/r_i})$
and that at the other is of form
$(t^{c_i/r_i}, 0, t^{\overline{ac_i}/r_i})$.\label{itm:one.one.IIb^2}
\item If $f$ is of type III
\textup{(\ref{txt:normIIb^2III}(\ref{itm:normIII}))},
then we may choose the coordinates so that
the restriction to $C^\dag$ at $Q$ is
of form $(t^{c/r}, 0, t^{\overline{ac}/r})$.\label{itm:one.one.III}
\end{enumerate}
\end{Lemma}

\begin{proof}
(\ref{itm:one.one.IIb^2})
If they both are of form $(0, t^{c_i/r_i}, t^{1-\overline{ac_i}/r_i})$,
then the image of the map
$\cO_Y(-aE) \otimes \cO_Y(E) \otimes \cO_Y((a-1)E) \otimes \cO_C \to \cO_C$
is $\fm_{\{Q_1, Q_2\} \subset C}$, whence
$s_C(-a)+s_C(1)+s_C(a-1)=-2$. Thus $s_C(-a)=0$ (\ref{txt:s_C}(\ref{itm:sIIb})),
and it is done by (\ref{prp:tfae}).
If they both are of form $(t^{c_i/r_i}, 0, t^{\overline{ac_i}/r_i})$,
then the image of the map
$\cO_Y(-aE) \otimes \cO_Y(-E) \otimes \cO_Y((a+1)E) \otimes \cO_C \to \cO_C$
is $\fm_{\{Q_1, Q_2\} \subset C}$, whence
$s_C(-a)+s_C(-1)+s_C(a+1)=-2$. Thus $s_C(-a)=0$ by $s_C(a+1) \le -1$ and
(\ref{txt:s_C}(\ref{itm:sIIb})), and it is done by (\ref{prp:tfae}).

(\ref{itm:one.one.III})
If the restriction is of form $(0, t^{c/r}, t^{1-\overline{ac}/r})$,
then the image of the map
$\cO_Y(-aE) \otimes \cO_Y(E) \otimes \cO_Y((a-1)E) \otimes \cO_C \to \cO_C$
is $\fm_{Q \subset C}$, whence
$s_C(-a)+s_C(1)+s_C(a-1)=-1$.
Thus $s_C(-a)=1$ (\ref{txt:s_C}(\ref{itm:sIIbIII})),
and it is done by (\ref{prp:tfae}).
\end{proof}

\begin{Proposition}\label{prp:-1-1}
\begin{enumerate}
\item Assume that $f$ is of type ${\textrm{IIb}}^{\vee \vee}$
with $s_C(-a)=-1$
\textup{(\ref{txt:normIIb^2III}(\ref{itm:normIIb^2}))}.
Then $a < r_1 < r_2$, $a_{11}=w_{Q_1}^C(1)=1$, $a_{22}=w_{Q_2}^C(-1)=1$
and $(H \cdot C)=1/r_1-1/r_2$.
In particular for each $i$ we may choose semi-invariant local coordinates
$x_{i1}, x_{i2}, x_{i3}$ with weights $\wt(x_{i1}, x_{i2}, x_{i3})=(1, -1, a)$
of the index one cover ${Q_i}^\sharp \in Y^\sharp$ so that
\begin{align*}
{Q_i}^\sharp \in C^\sharp \subset Y^\sharp
\cong o \in (x_{ii}\textrm{-}\mathrm{axis})
\subset x_{i1}x_{i2}x_{i3}\textrm{-}\mathrm{space}.
\end{align*}\label{itm:-1-1IIb^2}
\item Assume that $f$ is of type III
with $s_C(-a)=0$
\textup{(\ref{txt:normIIb^2III}(\ref{itm:normIII}))}.
Moreover we assume that $Q \in C$ has the local description in
\textup{(\ref{prp:gen.by.two}(\ref{itm:gen.two.ord}))}.
Then $a_1=w_{Q_1}^C(1)=1$ and $(H \cdot C)=1/r$.
In particular we may choose semi-invariant local coordinates
$x_1, x_2, x_3$ with weights $\wt(x_1, x_2, x_3)=(1, -1, a)$
of the index one cover $Q^\sharp \in Y^\sharp$ so that
\begin{align*}
Q^\sharp \in C^\sharp \subset Y^\sharp
\cong o \in (x_1\textrm{-}\mathrm{axis})
\subset x_1x_2x_3\textrm{-}\mathrm{space}.
\end{align*}\label{itm:-1-1III}
\end{enumerate}
\end{Proposition}

\begin{proof}
(\ref{itm:-1-1IIb^2})
First of all by (\ref{lem:one.one}(\ref{itm:one.one.IIb^2})),
we can choose $(A, B)=(1, 2)$ or $(2, 1)$ so that
the description at $Q_A$ is of form
$(t^{c_A/r_A}, 0, t^{\overline{ac_A}/r_A})$
and that at $Q_B$ is of form $(0, t^{c_B/r_B}, t^{1-\overline{ac_B}/r_B})$.

From the map
$\cO_Y(-E)^{\otimes r_Ar_B} \otimes \cO_Y(r_Ar_BE) \otimes \cO_C \to \cO_C$,
we obtain $r_Ar_Bs_C(-1)+r_Ar_B(E \cdot C)=-r_Bw_{Q_A}^C(1)-r_Aw_{Q_B}^C(1)$.
By $s_C(-1)=-1$ (\ref{txt:s_C}(\ref{itm:sIIb})),
$w_{Q_A}^C(1)=c_A$ and $w_{Q_B}^C(1)=r_B-c_B$, we have
\begin{align}\label{eqn:1+HC}
1+(H \cdot C)=\frac{c_A}{r_A}+\frac{r_B-c_B}{r_B}.
\end{align}
On the other hand from the map
$\cO_Y(-aE)^{\otimes r_Ar_B} \otimes \cO_Y(ar_Ar_BE) \otimes \cO_C \to \cO_C$,
we obtain $r_Ar_Bs_C(-a)+ar_Ar_B(E \cdot C)=-r_Bw_{Q_A}^C(a)-r_Aw_{Q_B}^C(a)$.
By $s_C(-a)=-1$,
$w_{Q_A}^C(a)=\overline{ac_A}$ and $w_{Q_B}^C(a)=r_B-\overline{ac_B}$, we have
\begin{align}\label{eqn:1+aHC}
1+a(H \cdot C)=\frac{\overline{ac_A}}{r_A}+\frac{r_B-\overline{ac_B}}{r_B}.
\end{align}
From these equalities (\ref{eqn:1+HC}) and (\ref{eqn:1+aHC}), one of the
inequalities $c_A < \overline{ac_A}$ and $c_B > \overline{ac_B}$ holds.
However $c_B > \overline{ac_B}$ cannot occur because
$(r_B, 0)$ is contained in the semigroup (\ref{eqn:semigroup}) for $Q_B$.
Hence we have $c_A < \overline{ac_A}$.

We claim that $c_A=1$. Suppose that $c_A \ge 2$
and take the smallest integer $m$ such that $c_Am \ge r_A$.
$w_{Q_A}^C(m) \equiv c_Am-r_A$ modulo $r_A$ and
$c_Am-r_A < a_{A1} < a_{A3}$. Hence the map
$\cO_Y(-mE) \otimes \cO_Y(mE) \otimes \cO_C \to \cO_C$ at $Q_A$ is
properly contained in $\fm_{Q_A \subset C}$.
In our situation the image of the map
$\cO_Y(aE) \otimes \cO_Y(-aE) \otimes \cO_C \to \cO_C$
is $\fm_{\{Q_1, Q_2\} \subset C}$.
Considering (\ref{txt:s_C}(\ref{itm:sIIb})) also, we have
\begin{enumerate}
\item $m>a$ or,\label{itm:m>a}
\item $r_1 \mid m$ and $m < a$.\label{itm:m<a}
\end{enumerate}
When (\ref{itm:m>a}) occurs, $\overline{ac_A}=ac_A$. Then $c_A=1$
since $(r_A, 0)$ is contained in the semigroup (\ref{eqn:semigroup}) for $Q_A$.
When (\ref{itm:m<a}) occurs, unless $m=2$ the map
$\cO_Y(-(m-1)E) \otimes \cO_Y((m-1)E) \otimes \cO_C \to \cO_C$ at $Q_A$ is
also properly contained in $\fm_{Q_A \subset C}$, which is a contradiction.
Thus the case where $m=2$ in (\ref{itm:m<a}) remains.
In this case $r_1=2$ and $c_1=1$.
Then we have $s_C(-a)=0$ as in the proof of (\ref{lem:one.one}),
which contradicts the assumption $s_C(-a)=-1$.

Now we have $(A, B)=(1, 2)$ and $r_1 < r_2$ by substituting $c_A=1$ in
(\ref{eqn:1+HC}). We have
\begin{align}
(H \cdot C)&=\frac{1}{r_1}-\frac{c_2}{r_2},\label{eqn:HC} \\
a(H \cdot C)&=\frac{\overline{a}}{r_1}-\frac{\overline{ac_2}}{r_2}.
\label{eqn:aHC}
\end{align}
Especially $\lfloor a/r_1 \rfloor = \lfloor ac_2/r_2 \rfloor$.

We claim that $c_2=1$ and hence $a < r_1$ (\ref{rem:a<r}) and
$(H \cdot C)=1/r_1-1/r_2$ (\ref{eqn:HC}).
We can show this similarly as we showed $c_1=1$ once we obtain
$c_2 < r_2-\overline{ac_2}$. But this inequality comes from 
the following obtained by (\ref{eqn:HC}) and (\ref{eqn:aHC}) respectively.
\begin{align*}
1-(H \cdot C)&=\frac{r_1-1}{r_1}+\frac{c_2}{r_2},\\
1+a(H \cdot C)&=\frac{\overline{a}}{r_1}+\frac{r_2-\overline{ac_2}}{r_2}.
\end{align*}

(\ref{itm:-1-1III})
By (\ref{lem:one.one}(\ref{itm:one.one.III}))
the description at $Q$ is of form
$(t^{c/r}, 0, t^{\overline{ac}/r})$ (\ref{prp:gen.by.two}).
From the map
$\cO_Y(-E)^{\otimes r} \otimes \cO_Y(rE) \otimes \cO_C \to \cO_C$,
we obtain $rs_C(-1)+r(E \cdot C)=-w_Q^C(1)$.
By $s_C(-1)=0$ (\ref{txt:s_C}(\ref{itm:sIIbIII})) and $w_Q^C(1)=c$, we have
$(H \cdot C)=c/r$.
On the other hand from the map
$\cO_Y(-aE)^{\otimes r} \otimes \cO_Y(arE) \otimes \cO_C \to \cO_C$,
we obtain $rs_C(-a)+ar(E \cdot C)=-w_Q^C(a)$.
By $s_C(-a)=0$ and $w_Q^C(a)=\overline{ac}$, we have
$a(H \cdot C)=\overline{ac}/r$.
Hence we have $ac = \overline{ac}$,
and we have the result because
$(r, 0)$ is contained in the semigroup (\ref{eqn:semigroup}).
\end{proof}

The next lemma is an easy consequence by numerical argument.

\begin{Lemma}\label{lem:irredp}
Assume that $f$ is of type ${\textrm{IIb}}^{\vee \vee}$.
Then $(H \cap E)_{\red}$ is irreducible.
\end{Lemma}

\begin{proof}
We may assume that $s_C(-a)=-1$ by (\ref{prp:tfae}).
Write $[H \cap E]=u[C]+[F]$ cycle theoretically,
where the support of $[F]$ does not pass $Q_1$ (\ref{txt:D}(\ref{itm:DIIb^2})).
Then we can write $(H \cdot [H \cap E])=u(H \cdot C)+v/r_2$.
By $(H \cdot [H \cap E])=E^3=a^{-1}(1/r_1+1/r_2)$ (\ref{thm:num.classn})
and $(H \cdot C)=1/r_1-1/r_2$ (\ref{prp:-1-1}(\ref{itm:-1-1IIb^2})), we obtain
\begin{align*}
\frac{au-1}{r_1}=\frac{au+1-av}{r_2}.
\end{align*}
Since $r_1 < r_2$ (\ref{prp:-1-1}(\ref{itm:-1-1IIb^2}))
and $a \ge 2$, we have $v=0$.
\end{proof}

\begin{Text}\label{txt:suitable}
We can take suitable surfaces by (\ref{prp:store}(\ref{itm:storeX2})).
\begin{enumerate}
\item Assume that $f$ is of type ${\textrm{IIb}}^{\vee \vee}$.
Then by (\ref{rem:di}(\ref{itm:diIIbIII})) and
(\ref{prp:store}(\ref{itm:storeX2}))
there exist surfaces $S$, $S_0$ on $Y$ defined by elements in $\cO_Y(-aE)$
such that $(S \cap E)_{\red} \not\subseteq
(H \cap E)_{\red} \cup (S_0 \cap E)_{\red}$.
Let $D$ be an irreducible reduced curve in $S \cap E$
which intersects $H$ and $S_0$ properly.
We can see that $D \cong \bP^1$ as in (\ref{lem:tree}).\label{itm:suiIIb^2}
\item Assume that $f$ is of type III.
Then by (\ref{rem:di}(\ref{itm:diIIbIII})) and
(\ref{prp:store}(\ref{itm:storeX2})) we can take a suitable $H$ and
a suitable $C \cong \bP^1 \subseteq H \cap E$ so that
$C$ properly intersects the strict transform $H_0$
of a general hyperplane section on $P \in X$.\label{itm:suiIII}
\end{enumerate}
\end{Text}

\begin{Lemma}\label{lem:irred}
Assume that $f$ is of type ${\textrm{IIb}}^{\vee \vee}$ with $s_C(-a)=-1$,
and take $C \subseteq H \cap E$ and $D \subseteq S \cap E$
as in \textup{(\ref{txt:normIIb^2III}(\ref{itm:normIIb^2}))}
and \textup{(\ref{txt:suitable}(\ref{itm:suiIIb^2}))}.
Then we can write $[S \cap E]=x[C]+[D]$ cycle theoretically
with a positive integer $x$, and
$D$ intersects $C$ exactly at one of $Q_1$ and $Q_2$.
Let $Q_i$ \textup{(}$i=1$ or $2$\textup{)} be the one $D$ passes through.
Then $(H \cdot D)=1/r_i$ and
we may choose semi-invariant local coordinates
$x_1, x_2, x_3$ with weights $\wt(x_1, x_2, x_3)=(1, -1, a)$
of the index one cover ${Q_i}^\sharp \in Y^\sharp$ so that
\begin{align*}
{Q_i}^\sharp \in D^\sharp \subset Y^\sharp
\cong o \in (x_1\textrm{-}\mathrm{axis})
\subset x_1x_2x_3\textrm{-}\mathrm{space}.
\end{align*}
\end{Lemma}

\begin{proof}
Write $[S \cap E]=x[C]+[D']$ cycle theoretically.
We can write $1/r_1+1/r_2=(H \cdot [S \cap E])=x(H \cdot C)+y/r_1+z/r_2$.
If $y=z=1$ then $[S \cap E]=[D']$ and $(S \cdot C) \ge 1/r_1+1/r_2$
because $S$ passes through $Q_1$ and $Q_2$.
On the other hand $(S \cdot C) \le (S \cdot [H \cap E])= 1/r_1+1/r_2$
(\ref{thm:num.classn}).
Hence $C=H \cap E$ and $s_C(-a)\ge0$ by
(\ref{rem:di}(\ref{itm:diIIbIII})) and (\ref{prp:store}(\ref{itm:storeX1X2}))
with $i=a<r_1$ (\ref{prp:-1-1}(\ref{itm:-1-1IIb^2})).
It contradicts $s_C(-a)=-1$.

Hence $y=0$ or $(y, z)=(1, 0)$.
By $(H \cdot [S \cap E])=1/r_1+1/r_2$ 
and $(H \cdot C)=1/r_1-1/r_2$ (\ref{prp:-1-1}(\ref{itm:-1-1IIb^2})),
we obtain
\begin{align*}
\frac{x-1+y}{r_1}=\frac{x+1-z}{r_2}.
\end{align*}
Since $y=0$ or $(y, z)=(1, 0)$, and $r_1 < r_2$
(\ref{prp:-1-1}(\ref{itm:-1-1IIb^2})), we have $(y, z)=(1, 0)$ or $(0, 1)$.
Therefore, provided that $r_1 \nmid r_2$, then
$D'=D$ and $D$ intersects $C$ only at $Q_1$ or $Q_2$
with $(H \cdot D)=1/r_1, 1/r_2$ respectively.
But if $r_1 \mid r_2$ and $(H \cdot D')=1/r_1$,
then setting $(r_1, r_2)=(r, sr)$,
by (\ref{lem:irredp}) we can write
\begin{align*}
au(\frac{1}{r}-\frac{1}{sr})=\frac{1}{r}+\frac{1}{sr},
\end{align*}
where $[H \cap E]=u[C]$ ($u \ge 2$).
By this equation and
$s > 1$ (\ref{prp:-1-1}(\ref{itm:-1-1IIb^2})),
we have $(s, au)=(2, 3)$ or $(3, 2)$.
But this contradicts $a \ge 2$ and $u \ge 2$.

Let $Q_i$ ($i=1$ or $2$) be the one $D$ passes through.
From the map
$\cO_Y(-E)^{\otimes r_i} \otimes \cO_Y(r_iE) \otimes \cO_D \to \cO_D$,
we obtain $r_is_D(-1)+r_i(E \cdot D)=-w_{Q_i}^D(1)$.
Hence $w_{Q_i}^D(1) \equiv 1$ modulo $r_i$.
Therefore $w_{Q_i}^D(1) = 1$ since
$(r_i, 0)$ is contained in the semigroup (\ref{eqn:semigroup}) for $Q_i \in D$,
whence we have the desired description.
\end{proof}

The next theorem gives complete description in
types ${\textrm{IIb}}^{\vee \vee}$ and III.

\begin{Theorem}\label{thm:localIIb^2III}
Assume that $f$ is of type ${\textrm{IIb}}^{\vee \vee}$ and III.
Then there exist an $H$ and a curve $C \subseteq H \cap E$ in
\textup{(\ref{txt:normIIb^2III})}
which satisfy the assumption in \textup{(\ref{prp:tfae})}.
\end{Theorem}

\begin{proof}
We first treat type ${\textrm{IIb}}^{\vee \vee}$.
Suppose that $C$ does not satisfy the assumption in (\ref{prp:tfae}), that is,
suppose that $s_C(-a)=-1$ (\ref{txt:s_C}(\ref{itm:sIIb})).
Then by (\ref{txt:suitable}(\ref{itm:suiIIb^2})) and (\ref{lem:irred}),
there exist an $S$ and a $Q_i \in D \cong \bP^1$
in (\ref{txt:suitable}(\ref{itm:suiIIb^2}))
such that we may choose semi-invariant local coordinates $x_1, x_2, x_3$
so that $(x_1, x_2, x_3) |_{D^\dag}=(t^{1/r_i}, 0, 0)$.
$D$ intersects $H$ properly.

Let $H^\sharp$ be the preimage of $H$ on the index one cover
$Q_i^\sharp \in Y^\sharp$.
Write the defining equation of $H^\sharp$ at $Q_i^\sharp$ as
$h(x_1)+g_1x_2+g_2x_3$,
where $h(x_1)$ is a semi-invariant formal series of weight $1$
and $g_1$, $g_2$ are semi-invariant formal series of $x_1, x_2, x_3$.
The order of the equation of $H^\sharp |_{D^\dag}$ with respect to $t$
equals $(H \cdot D)=1/r_i$ (\ref{lem:irred}).
It means that the order of $h(t^{1/r_i})$ with respect to $t$ is $1/r_i$,
whence the coefficient of $x_1$ in $h$ is not zero.
Hence we may choose the coordinates $x_1, x_2, x_3$ so that
$H^\sharp$ is given by $x_1=0$.
If $Q_i=Q_1$ then for the coordinates $x_{11}, x_{12}, x_{13}$
in (\ref{prp:-1-1}(\ref{itm:-1-1IIb^2}))
we can write $x_1=cx_{11}+g(x_{11}, x_{12}, x_{13})$,
where $c$ is not zero and $g$ is a semi-invariant formal series of weight $1$
which does not contain the term $x_{11}$.
Thus we can easily to see that $x_1 |_{C^\dag} \neq 0$
by $(x_{11}, x_{12}, x_{13}) |_{C^\dag}=(t^{1/r_1}, 0, 0)$,
which contradicts $C \subset H$.

Hence $Q_i=Q_2$.
We consider a general $S_0$ in (\ref{txt:suitable}(\ref{itm:suiIIb^2})).
$S_0 \cap E$ also has a component $D_0 \cong \bP^1$ which has the
same properties as $D \subset S \cap E$ has, and we can choose new coordinates
${x_1}', {x_2}', {x_3}'$ and $t'$ for $Q_2 \in D_0 \subset Y$
so that $({x_1}', {x_2}', {x_3}') |_{D_0^\dag}=({t'}^{1/r_2}, 0, 0)$.
Since $(S \cdot D_0)=a(H \cdot D_0)=a/r_2<1$ (\ref{lem:irred})(\ref{rem:a<r}),
$S$ intersects $D_0$ only at $Q_2$.
Let $S^\sharp$ be the preimage of $S$ on $Y^\sharp$,
and let $g(x_1, x_2, x_3)=h({x_1}', {x_2}', {x_3}')$ be
the defining equation of $S^\sharp$ at $Q_2^\sharp$.
The order of $h({t'}^{1/r_2}, 0, 0)$, the equation of $S^\sharp |_{D_0^\dag}$,
with respect to $t$
equals $(S \cdot D_0)=a/r_2$, whence we can write
$h=(c'+h_1({x_1}'){x_1'}^{r_2}){x_1'}^a+h_2({x_1}', {x_2}', {x_3}'){x_2}'
+h_3({x_1}', {x_2}', {x_3}'){x_3}'$,
where $c' \in \bC^\times$ and $h_1, h_2, h_3$ are semi-invariant
with weights $0, a+1, 0$.
Thus we can write $g=cx_1^a+ \cdots$ ($c \in \bC^\times$)
and we see that $g(t^{1/r_2}, 0, 0) = 0$ holds
only if $h_3$ is a unit, by considering weights
and $a<r_1<r_2$ (\ref{txt:suitable}(\ref{itm:suiIIb^2})).
Of course $g(t^{1/r_2}, 0, 0) = 0$ by $D \subset S$, whence
$h_3$ is a unit and $S^\sharp$ is smooth at $Q_2^\sharp$.

Therefore we may choose the coordinates
$x_{21}, x_{22}, x_{23}$ in (\ref{prp:-1-1}(\ref{itm:-1-1IIb^2}))
so that $x_{21}=x_1$, the equation of $H^\sharp$, and that
$x_{23}$ is the equation of $S^\sharp$. Especially the scheme
$S \cap H$ is irreducible and reduced at the generic point of $C$.
Write $[H \cap E]=u[C]$ (\ref{lem:irredp}) and
$[S \cap E]=x[C]+[D]$ (\ref{lem:irred}).
Then  we have
\begin{align*}
ua \Bigl( \frac{1}{r_1}-\frac{1}{r_2} \Bigr) &= \frac{1}{r_1}+\frac{1}{r_2}, \\
x \Bigl( \frac{1}{r_1}-\frac{1}{r_2} \Bigr) + \frac{1}{r_2}
&= \frac{1}{r_1}+\frac{1}{r_2}.
\end{align*}
Thus $(ua-x)(r_2-r_1)=r_1$ and we can write $r_1=nr$, $r_2=(n+1)r$.
We note that
\begin{align}
ua&=2n+1,\label{eqn:ua=} \\
x&=n+1.\label{eqn:x=}
\end{align}

On the other hand let $\pi \colon \hat{E} \to E$ be the normalisation of $E$
and set $\hat{C}=(C \times_E \hat{E})_{\red}$.
The coefficients of $[C]$ in the 1-cycles $[S \cap E]$ and $[aH \cap E]$
are determined by the lengths of the schemes
$\pi^*S \cap \hat{C}$ and $\pi^*(aH) \cap \hat{C}$ at the generic points
of all the irreducible components of $\hat{C}$.
From this point of view, we see that $(S \cap E) \cap (aH \cap E)$ contains
$1$-cycle $x[C]$ because of the general choice of $S$.
$(S \cap E) \cap (aH \cap E) \subseteq S \cap aH$ and
$S \cap H$ is irreducible and reduced at the generic point of $C$.
Hence $x \le a$ but it contradicts
(\ref{eqn:ua=}), (\ref{eqn:x=}) and $u \ge 2$.

We now treat type III.
We start with $H$, $Q \in C$ and $H_0$
in (\ref{txt:suitable}(\ref{itm:suiIII})).
Suppose that $C$ does not satisfy the assumption in (\ref{prp:tfae}),
that is, suppose that $s_C(-a)=0$.
Then by (\ref{prp:gen.by.two}) and (\ref{prp:-1-1}(\ref{itm:-1-1III}))
we may choose coordinates $x_1, x_2, x_3$
so that $(x_1, x_2, x_3) |_{C^\dag}=(t^{1/r}, 0, 0)$,
$(t^{2/r}, t^{1-2/r}, 0)$ or $(t^{1-2/r}, t^{2/r}, 0)$.
We note that $(H_0 \cdot E)=1/r, 2/r, 1-2/r$ respectively in these cases.
Thus $H_0$ intersects $C$ only at $Q$.

Let $H_0^\sharp$ be the preimage of $H_0$ on the index one cover
$Q^\sharp \in Y^\sharp$ and let $g(x_1, x_2, x_3)$
be the defining equation of $H_0^\sharp$ at $Q^\sharp$.
$g$ is semi-invariant of weight $1$ and the order of $g(t^{1/r}, 0, 0)$,
$g(t^{2/r}, t^{1-2/r}, 0)$ or $g(t^{1-2/r}, t^{2/r}, 0)$ with respect to $t$
is $1/r, 2/r, 1-2/r$ respectively.
As in the former argument we can see that
the coefficient of $x_1$ in the equation of $H_0^\sharp$
is not zero in each of the cases.
Hence we may choose semi-invariant local coordinates
$x_1, x_2, x_3$ with weights $\wt(x_1, x_2, x_3)=(1, -1, a)$
of $Q^\sharp \in Y^\sharp$ so that $H_0$ is given by $x_1=0$.

Now we take $Q \in C_0 \cong \bP^1 \subset H_0 \cap E$.
If $s_{C_0}(-a)\neq 1$, then $s_{C_0}(-a)=0$ and again
we may choose new coordinates ${x_1}', {x_2}', {x_3}'$ and $t'$
so that $({x_1}', {x_2}', {x_3}') |_{C_0^\dag}=({t'}^{1/r}, 0, 0)$,
$({t'}^{2/r}, {t'}^{1-2/r}, 0)$ or $({t'}^{1-2/r}, {t'}^{2/r}, 0)$.
On the other hand we can write $x_1=c{x_1}'+g({x_1}', {x_2}', {x_3}')$,
where $c$ is not zero and $g$ is a semi-invariant formal series of weight $1$
which does not contain the term $x_1'$.
Thus we can easily see that $x_1 |_{C_0^\dag} \neq 0$,
which contradicts $C_0 \subset H_0$.
\end{proof}

\section{Existence of Du Val sections}\label{sec:existence}
\begin{Text}\label{txt:review}
In this section we show the existence of Du Val sections in the
anticanonical system of $Y$ (\ref{thm:gen.element})
together with providing more information.
It is done when $f$ is of type O or I (\ref{txt:easy.elephant}),
and we have nothing to do in type IV (\ref{rem:typeIV}).
Thus we keep the assumption that $f$ is of type II or III.
First we restate results in Section \ref{sec:local}.
\begin{enumerate}
\item (type IIa (\ref{thm:localIIa}))
\begin{enumerate}
\item $a=2$ and $[H \cap E]=[\bP^1]+[\bP^1]$.
$s_C(-2)=0$ and $(H \cdot C)=1/r$.
$Q \in C$ is locally expressed by semi-invariant local coordinates
$x_1, x_2, x_3$ with weights $\wt(x_1, x_2, x_3)=(1, -1, 4)$
of the index one cover $Q^\sharp \in  Y^\sharp$ so that
$(x_1, x_2, x_3) |_{C^\dag}=(t^{(r+1)/2r}, 0, t^{2/r})$
($r \equiv 1$ modulo $4$) or $(0, t^{(r-1)/2r}, t^{2/r})$
($r \equiv 3$ modulo $4$).
\item $r=5$, $a=4$ and $H \cap E=\bP^1$. $s_C(-4)=0$ and $(H \cdot C)=1/r$.
$Q \in C$ is locally expressed by semi-invariant local coordinates
$x_1, x_2, x_3$ with weights $\wt(x_1, x_2, x_3)=(1, -1, 3)$
of the index one cover $Q^\sharp \in  Y^\sharp$
as $(x_1, x_2, x_3) |_{C^\dag}=(t^{3/5},t^{2/5}, 0)$.
\end{enumerate}\label{itm:revIIa}
\item (type ${\textrm{IIb}}^{\vee}$ (\ref{thm:localIIb^}))
\begin{enumerate}
\item $J=\{(r, 1), (r, 1)\}$, $a=2$ and $H \cap E=\bP^1$.
$s_C(-2)=0$ and $(H \cdot C)=1/r$.
$Q \in C$ is locally expressed by semi-invariant local coordinates
$x_1, x_2, x_3, x_4$ with weights $\wt(x_1, x_2, x_3, x_4)=(1, -1, 2, 0)$
of the index one cover $Q^\sharp \in Y^\sharp$
as $(x_1, x_2, x_3, x_4) |_{C^\dag}=(0, 0, t^{2/r}, t)$.
\item $J=\{(2, 1), (4, 1)\}$, $a=3$ and $H \cap E=\bP^1$.
$s_C(-3)=0$ and $(H \cdot C)=1/4$.
$Q \in C$ is locally expressed by semi-invariant local coordinates
$x_1, x_2, x_3, x_4$ with weights $\wt(x_1, x_2, x_3, x_4)=(1, 3, 3, 2)$
of the index one cover $Q^\sharp \in Y^\sharp$
as $(x_1, x_2, x_3, x_4) |_{C^\dag}=(0, 0, t^{3/4}, t^{1/2})$.
\end{enumerate}\label{itm:revIIb^}
\item (type ${\textrm{IIb}}^{\vee \vee}$ (\ref{thm:localIIb^2III}))
$H \cap E=\bP^1$.
$s_C(-a)=0$ and $(H \cdot C)=(r_1+r_2)/ar_1r_2$.
$Q_i \in C$ is locally expressed by
semi-invariant local coordinates $x_{i1}, x_{i2}, x_{i3}$ with weights
$\wt(x_{i1}, x_{i2}, x_{i3})=(1, -1, a)$
of the index one cover $Q_i^\sharp \in Y^\sharp$
as $(x_{i1}, x_{i2}, x_{i3}) |_{C^\dag}=(0, 0, t^{1/r_i})$.\label{itm:revIIb^2}
\item (type III (\ref{thm:localIIb^2III}))
We can choose $H$ so that $H \cap E=\bP^1$.
$s_C(-a)=1$ and $(H \cdot C)=(r+1)/ar$.
$Q \in C$ is locally expressed by
semi-invariant local coordinates $x_1, x_2, x_3$ with weights
$\wt(x_1, x_2, x_3)=(1, -1, a)$
of the index one cover $Q^\sharp \in Y^\sharp$
as $(x_1, x_2, x_3) |_{C^\dag}=(0, 0, t^{1/r})$.\label{itm:revIII}
\end{enumerate}
\end{Text}

\begin{Text}\label{txt:S}
Let $S$ be a surface on $Y$ defined by a general element in $\cO_Y(-aE)$
and let $S_X$ be its strict transform on $X$.
Our main goal is to see that $S$ has at worst Du Val singularities,
but we furthermore induce some information on $S \cap E$.
\end{Text}

\begin{Lemma}\label{lem:ex}
Assume that $f$ is of type IIa with $a=2$, and
write $[H \cap E]=[C]+[C']$ cycle theoretically.
Then there exists a natural exact sequence
\begin{align*}
0 \to \cO_{C'}(-1) \to \cO_{H \cap E} \to \cO_C \to 0.
\end{align*}
\end{Lemma}

\begin{proof}
Let $\cI$ be the kernel of the natural map
$\cO_{H \cap E} \twoheadrightarrow \cO_C$.
$\cI$ is supported on $C' \cong \bP^1$.
By $H^0(\cO_{H \cap E})=\bC$ (\ref{rem:tree}), we have $H^0(\cI)=0$.
By $H^1(\cO_{H \cap E})=0$ (\ref{lem:tree}(\ref{itm:curvevan}))
we have $H^1(\cI)=0$, whence $\cI \cong \cO_{\bP^1}(-1)$.
\end{proof}

\begin{Theorem}\label{thm:el.irred}
\begin{enumerate}
\item If $f$ is of type II, then $S$ intersects $H \cap E$
at and only at the non-Gorenstein points of $Y$.\label{itm:el.II}
\item If $f$ is of type III, then for any $Q' \in C=H \cap E$
which is a smooth point of $Y$, there exists an $S$
such that $S$ intersects $H \cap E$
at and only at $Q'$ and the non-Gorenstein point of $Y$.\label{itm:el.III}
\end{enumerate}
\end{Theorem}

\begin{proof}
(\ref{itm:el.II})
Take any $C \cong \bP^1 \subseteq H \cap E$.
By (\ref{eqn:vanish}), (\ref{txt:comm.diag}) and (\ref{lem:ex})
we obtain a surjective map
$f_*\cO_Y(-aE) \twoheadrightarrow  H^0([\cO_Y(-aE)]_C)$.
Recalling that $s_C(-a)=0$ (\ref{txt:review}),
we have $H^0([\cO_Y(-aE)]_C) = \bC$.
Hence we find an $S$ which properly intersects $H \cap E$.
Of course $S$ passes through any non-Gorenstein point of $Y$.
On the other hand because
$(S \cdot [H \cap E])=aE^3 \le 1$ in any case (\ref{txt:review}),
we see that $S$ does not intersect $H \cap E$ at any other point of $Y$.

(\ref{itm:el.III})
It suffices to consider $[\cO_Y(-aE) \otimes \fm_{Q' \subset Y}]_C$
instead of $[\cO_Y(-aE)]_C$ since $s_C(-a)=1$ (\ref{txt:review}).
\end{proof}

It is the time to show our main theorem (\ref{thm:gen.element}).

\begin{Theorem}\label{thm:II.III.elephant}
Assume that $f$ is of type II or III.
Then for a suitable $S$ in \textup{(\ref{thm:el.irred})},
$S$ has at worst Du Val singularities.
Moreover the type of a Du Val singularity at any $Q \in S$ is that of
a general Du Val section of a germ $Q \in Y$.
\end{Theorem}

\begin{proof}
We have only to consider the case where $P$ is a singular point
by \cite{Ka01}. Because of (\ref{thm:el.irred}) and Bertini's theorem,
it is enough to show that
$S$ has a Du Val singularity at each non-Gorenstein point $Q \in Y$
whose type equals that of a general Du Val section of a germ $Q \in Y$.

First we treat type IIa (\ref{txt:review}(\ref{itm:revIIa})).
Let $Q \in C = H \cap E$ be the non-Gorenstein point with a curve
investigated.

Consider the case $a=2$.
It is enough to show that the coefficient of $x_3$ in
the semi-invariant equation $h(x_1, x_2, x_3)$ of $S$ in (\ref{txt:S})
is not zero.
Since $S$ intersects $C$ only at $Q$ (\ref{thm:el.irred}),
the order of $h(t^{(r+1)/2r}, 0, t^{2/r})$
($r \equiv 1$ modulo $4$) or $h(0, t^{(r-1)/2r}, t^{2/r})$
($r \equiv 3$ modulo $4$) with respect to $t$ equals
$(S \cdot C)=2/r$.
Hence the coefficient of $x_3$ in $h$ is not zero.

Consider the case $a=4$.
$-K_Y$ is linearly equivalent to $E+5H$. Thus it is enough to show that
\begin{enumerate}
\item $E$ has a Du Val singularity of type $A_4$ at $Q$, and\label{itm:A4}
\item $\left| 5H \right|$ is free at $Q$.\label{itm:5H}
\end{enumerate}

To see (\ref{itm:A4}), consider a birational morphism $g \colon Z \to Y$
such that $Z$ has a $g$-exceptional divisor $F$
whose discrepancy with respect to $K_X$ is $1$ \cite{Ma96}. Write
\begin{align*}
K_Z&=g^*K_Y+bF+(\textrm{others}),\\
g^*E&=E_Z+mF+(\textrm{others}),
\end{align*}
where $E_Z$ is the strict transform of $E$.
Then $K_Z=g^*(f^*K_X+4E)+bF+(\textrm{others})=g^*f^*K_X+4E_Z+
(b+4m)F+(\textrm{others})$, whence $b=m=1/5$ by $5b, 5m \in \bZ$.
Hence the defining equation of $E$ by $x_1, x_2, x_3$
in (\ref{txt:review}(\ref{itm:revIIa})) has a nonzero linear term.
Considering the weight we see that the coefficient of $x_3$ in
the equation of $E$ is not zero, whence we may take $x_3$ so that $E$
is given by $x_3=0$. This implies (\ref{itm:A4}).

It is easy to see (\ref{itm:5H}).
From (\ref{txt:comm.diag}) and (\ref{eqn:vanish})
we obtain a surjective map
$f_*\cO_Y(-5E) \twoheadrightarrow H^0([\cO_Y(-5E)]_C)
=H^0(\cO_{\bP^1}(1))$, which shows (\ref{itm:5H}).

Secondly we treat type ${\textrm{IIb}}^{\vee}$
(\ref{txt:review}(\ref{itm:revIIb^})).
Let $Q \in C = H \cap E$ be the non-Gorenstein point with a curve
investigated.
$Q \in Y$ is described as in (\ref{rem:precise}).
Consider a surface $T$ on a germ $Q \in Y$
defined by a section in $\cO_Y(-K_Y)_Q$
whose equation in $y_1, y_2, y_3, y_4$ contains a linear term.
Then by the list (\ref{rem:precise}) we have
\begin{enumerate}
\item If $J=\{(r, 1), (r, 1)\}$ and $r \ge 5$, then
the linear term of the equation of $T$ is $y_3$.
Any such $T$ gives a Du Val singularity of type $A_{2r-1}$,
which is the type of a general Du Val section of a germ $Q \in Y$.
\item In any other case, the linear term of the equation of $T$ is
a linear combination of $y_2$ or $y_3$.
Any $T$ equipped with a general linear term
gives a Du Val singularity of type $A_5$, $E_6$ or $D_5$
according to (\ref{itm:A/r}), (\ref{itm:D/3}) or (\ref{itm:A/4})
in (\ref{rem:precise}) respectively,
which is the type of a general Du Val section of a germ $Q \in Y$.
\end{enumerate}
Thus if $J=\{(r, 1), (r, 1)\}$ and $r \ge 5$, it is enough to show that
\begin{sText}\label{txt:eq.S}
The equation of $S$ in (\ref{txt:S}) by $y_1, y_2, y_3, y_4$
contains a linear term.
\end{sText}
In any other case, $-K_Y$ is linearly equivalent to $E+rH$.
Thus it is enough to show that, adding to (\ref{txt:eq.S}),
\begin{sText}\label{txt:to.show.add}
\begin{enumerate}
\item The equation of $E$ by $y_1, y_2, y_3, y_4$
contains a linear term which is linearly independent to that in
the equation of $S$.\label{itm:to.show.add.E}
\item $\left| rH \right|$ is free at $Q$.\label{itm:to.show.add.H}
\end{enumerate}
\end{sText}

However in (\ref{txt:to.show.add}),
the existence of a linear term in the equation of $E$ and
the statement (\ref{txt:to.show.add}(\ref{itm:to.show.add.H})) can be shown
similarly as in the proof of type IIa with $a=4$.
In particular in (\ref{txt:to.show.add})
we may choose coordinates $x_1, x_2, x_3, x_4$ in
(\ref{txt:review}(\ref{itm:revIIb^})) so that $E$ is given by $x_2=0$.
Hence (\ref{txt:eq.S}) and (\ref{txt:to.show.add}) follows if we prove
that the coefficient of $x_3$ in the defining equation
$h(x_1, x_2, x_3, x_4)$ of $S$ in (\ref{txt:S}) is not zero.

Because $S$ intersects $C$ only at $Q$ (\ref{thm:el.irred}),
the order of $h(0, 0, t^{2/r}, t)$ or $h(0, 0, t^{3/4}, t^{1/2})$
with respect to $t$ equals $(S \cdot C)=aE^3=2/r$ or $3/4$
respectively according to $J=\{(r, 1), (r, 1)\}$ or $\{(2, 1), (4, 1)\}$.
Hence the coefficient of $x_3$ in $h$ is not zero.

Finally we treat types ${\textrm{IIb}}^{\vee \vee}$ and III
(\ref{txt:review}(\ref{itm:revIIb^2})(\ref{itm:revIII})).
Let $Q \in C = H \cap E$ be a non-Gorenstein point with a curve
investigated. It is enough to show that the coefficient of $x_3$ in the
defining equation $h(x_1, x_2, x_3)$ of $S$ in (\ref{txt:S}) is not zero, and
we can see this similarly, remarking that the local intersection number
of $S$ and $C$ at $Q$ is $1/r_Q$
(\ref{txt:review}(\ref{itm:revIIb^2})(\ref{itm:revIII}))(\ref{thm:el.irred}).
\end{proof}

\section{Possible types of singularities}\label{sec:possible}
\begin{Text}\label{txt:graph}
In this section we restrict types of Du Val singularities
on $S$ and $S_X$, and prove (\ref{thm:detail}) as a consequence.
First we recall dual graphs for
minimal resolutions of Du Val singularities.
$\circ$ denotes an exceptional curve,
$\bullet$ denotes the strict transform of a general hyperplane section, and
numbers attached exceptional curves $F_i$
are the coefficients of $F_i$ in the fundamental cycle.
\begin{enumerate}
\item (type $A_n$)
\begin{align*}
\begin{array}{ccccccccc}
\bullet & - & \overset{F_1}{\underset{1}{\circ}} & - &
\cdots  & - & \overset{F_n}{\underset{1}{\circ}} & - & \bullet
\end{array}
\end{align*}
\item (type $D_n$)
\begin{align*}
\begin{array}{ccccccccc}
& & \bullet   & & & & \overset{F_n}{\underset{1}{\circ}} & & \\
& & \shortmid & & & & \shortmid                          & & \\
\overset{F_1}{\underset{1}{\circ}} & - & \overset{F_2}{\underset{2}{\circ}} &
- & \cdots  & - & \overset{F_{n-2}}{\underset{2}{\circ}} &
- & \overset{F_{n-1}}{\underset{1}{\circ}}
\end{array}
\end{align*}
\item (type $E_6$)
\begin{align*}
\begin{array}{ccccccccc}
& & & & \bullet                            & & & & \\
& & & & \shortmid                          & & & & \\
& & & & \overset{F_1}{\underset{2}{\circ}} & & & & \\
& & & & \shortmid                          & & & & \\
\overset{F_2}{\underset{1}{\circ}} & - & \overset{F_3}{\underset{2}{\circ}}
& - &
\overset{F_4}{\underset{3}{\circ}} & - & \overset{F_5}{\underset{2}{\circ}}
& - & \overset{F_6}{\underset{1}{\circ}}
\end{array}
\end{align*}
\item (type $E_7$)
\begin{align*}
\begin{array}{ccccccccccccc}
& & & & & & \overset{F_7}{\underset{2}{\circ}} & & & & & & \\
& & & & & & \shortmid                          & & & & & & \\
\bullet & - & \overset{F_1}{\underset{2}{\circ}} & - &
\overset{F_2}{\underset{3}{\circ}} & - & \overset{F_3}{\underset{4}{\circ}}
& - &
\overset{F_4}{\underset{3}{\circ}} & - & \overset{F_5}{\underset{2}{\circ}}
& - & \overset{F_6}{\underset{1}{\circ}}
\end{array}
\end{align*}
\item (type $E_8$)
\begin{align*}
\begin{array}{ccccccccccccccc}
& & & & & & & & & & \overset{F_8}{\underset{3}{\circ}} & & & & \\
& & & & & & & & & & \shortmid                          & & & & \\
\bullet & - & \overset{F_1}{\underset{2}{\circ}} & - &
\overset{F_2}{\underset{3}{\circ}} & - & \overset{F_3}{\underset{4}{\circ}}
& - &
\overset{F_4}{\underset{5}{\circ}} & - & \overset{F_5}{\underset{6}{\circ}}
& - &
\overset{F_6}{\underset{4}{\circ}} & - & \overset{F_7}{\underset{2}{\circ}}
\end{array}
\end{align*}
\end{enumerate}
\end{Text}

\begin{Text}\label{txt:part.resol}
Take a surface $S$ on $Y$ defined by a general element in $\left| -K_Y \right|$
(\ref{thm:O.I.elephant})(\ref{thm:II.III.elephant}), and
define $S_X$ as its strict transform on $X$.
A general hyperplane section $P \in H_X$ on $X$ gives
also a general hyperplane section
$P \in D_X=H_X |_{S_X}$ on $S_X$.
Let $f_S \colon S \to S_X$ be the induced map.
Since $S$ and $S_X $ have at worst Du Val singularities
(\ref{thm:O.I.elephant})(\ref{thm:II.III.elephant}),
$f_S$ factors the minimal resolution of $S_X$.
Let $Z$ be the strict transform on $S$ of the fundamental cycle.
$f^*H_X=H+bE$ in (\ref{txt:easy.elephant}) and
$f^*D_X=D+Z=H |_S + bE |_S$,
where $D$ is the strict transform on $S$ of $D_X$.
We note that $bE |_S \le Z$. We also have
\begin{align}\label{eqn:key.value}
(H |_S \cdot bE |_S)=b^2aE^3.
\end{align}
\end{Text}

\begin{Text}\label{txt:typeI}
Assume that $f$ is of type I. In this case one of the following holds
by (\ref{thm:num.classn}), (\ref{prp:formulae}) and (\ref{thm:O.I.elephant}).
\begin{enumerate}
\item $J=\{(7, 3)\}$, $E^3=1/7$ and $a=b=2$.\label{itm:typeI.7}
\item $J=\{(3, 1), (5, 2)\}$, $E^3=1/15$ and $a=b=2$.\label{itm:typeI.15}
\end{enumerate}
\end{Text}

\begin{Lemma}\label{lem:I.graph}
Assume that $f$ is of type I.
Take a general $S$ in \textup{(\ref{txt:part.resol})}.
\begin{enumerate}
\item If $J=\{(7, 3)\}$ \textup{(\ref{txt:typeI}(\ref{itm:typeI.7}))},
then the dual graph for the partial resolution $S \to S_X$ is
\begin{align*}
\ast - \underset{2}{\circ} \quad \textrm{or}
\quad \underset{2}{\circ} - \ast - \underset{2}{\circ},
\end{align*}
where the attached numbers are the coefficients in $2E |_S$,
and $\ast$ means that a Du Val singularity of type at best $A_6$
appears there.\label{itm:lem.typeI.7}
\item If $J=\{(3, 1), (5, 2)\}$ \textup{(\ref{txt:typeI}(\ref{itm:typeI.15}))},
then the dual graph for the partial resolution $S \to S_X$ is
\begin{align*}
\ast - \underset{2}{\circ} - \star \quad \textrm{or}
\quad \ast - \underset{4}{\circ} - \star,
\end{align*}
the attached number is the coefficient in $2E |_S$,
and $\ast$, $\star$ means that a Du Val singularity of type
at best $A_2$, $A_4$ appears there respectively.\label{itm:lem.typeI.15}
\end{enumerate}
\end{Lemma}

\begin{proof}
(\ref{itm:lem.typeI.7})
Since $(-E \cdot [S \cap E])=2E^3=2/7$,
it is enough to show that the existence of $S$ such that
$S \cap E$ defines a reduced $1$-cycle.
In this case we can calculate $\chi(\cQ_{-2})=2$ by (\ref{eqn:A_i}).
Hence we have a desired $S$ as in the proof of
(\ref{prp:store}(\ref{itm:storeX2})).

(\ref{itm:lem.typeI.15})
Since $(-E \cdot [S \cap E])=2E^3=2/15$,
it is enough to show that the existence of $S$ such that
$S \cap E$ defines an irreducible but possibly nonreduced $1$-cycle.
By $(-E \cdot [S \cap E])<1/5$, each irreducible component of
$(S \cap E)_{\red}$ passes through all the non-Gorenstein points of $Y$.
Hence $(S \cap E)_{\red}=\bP^1$ by $h^1(\cO_{S \cap E})=0$ (\ref{rem:tree}).
\end{proof}

\begin{Text}\label{txt:list.typeI}
We restrict the dual graph for the partial resolution $S$ in
(\ref{txt:part.resol}) considering $2E |_S \le Z$ and (\ref{lem:I.graph}).
The following list shows all the possibilities up to permutation
in terms of $F_i$ in (\ref{txt:graph}).
\begin{center}
\begin{tabular}[t]{l|l|l}
\multicolumn{3}{c}{$J=\{(7, 3)\}$ (\ref{txt:typeI}(\ref{itm:typeI.7}))} \\
\hline
\multicolumn{1}{c|}{case} & \multicolumn{1}{c|}{$S_X$}
& \multicolumn{1}{c}{$Z$} \\
\hline
(a) & $D_n$ & $Z \not\supseteq F_1,\ldots,F_{m-1}$ \\
    &       & and $Z \supseteq 2F_m$ \\
(b) & $E_7$ & $Z \supseteq 2F_1$ \\
(c) & $E_7$ & $Z=2F_5$ \\
(d) & $E_7$ & $Z=2F_7$ \\
(e) & $E_8$ & $Z \supseteq 2F_1$ \\
(f) & $E_8$ & $Z=3F_2$ \\
(g) & $E_8$ & $Z=3F_2+2F_7$ \\
(h) & $E_8$ & $Z=4F_3$ \\
(i) & $E_8$ & $Z=4F_3+2F_7$ \\
(j) & $E_8$ & $Z=4F_6$ \\
(k) & $E_8$ & $Z=2F_7$ \\
(l) & $E_8$ & $Z=2F_7+3F_8$ \\
(m) & $E_8$ & $Z=3F_8$ \\
\hline
\end{tabular} \qquad
\begin{tabular}[t]{l|l|l}
\multicolumn{3}{c}{$J=\{(3, 1), (5, 2)\}$
(\ref{txt:typeI}(\ref{itm:typeI.15}))} \\
\hline
\multicolumn{1}{c|}{case} & \multicolumn{1}{c|}{$S_X$}
& \multicolumn{1}{c}{$Z$} \\
\hline
(a) & $D_n$ & $Z \not\supseteq F_1,\ldots,F_{m-1}$ \\
    &       & and $Z \supseteq 2F_m$ \\
(b) & $E_7$ & $Z=3F_4$ \\
(c) & $E_8$ & $Z=4F_3$ \\
(d) & $E_8$ & $Z=5F_4$ \\
(e) & $E_8$ & $Z=6F_5$ \\
\hline
\end{tabular}
\end{center}

We calculate $(H |_S \cdot 2E |_S)$ for all the possible cases.
\begin{center}
\begin{tabular}[t]{l|l|l}
\multicolumn{3}{c}{$J=\{(7, 3)\}$} \\
\hline
\multicolumn{1}{c|}{case} &
\multicolumn{1}{c|}{$2E_S$} & \multicolumn{1}{c}{$(H |_S \cdot 2E |_S)$} \\
\hline
(a) & $2F_m+\cdots$ & $4/m$ \\
(b) & $2F_1+\cdots$ & $2$ \\
(c) & $2F_5$        & $1$ \\
(d) & $2F_7$        & $8/7$ \\
(e) & $2F_1+\cdots$ & $2$ \\
(f) & $2F_2$        & $2/3$ \\
(g) & $2F_2+2F_7$   & $1$ \\
(h) & $2F_3$        & $1/3$ \\
(i) & $2F_3+2F_7$   & $4/3$ \\
(j) & $2F_6$        & $2/7$ \\
(k) & $2F_7$        & $1$ \\
(l) & $2F_7+2F_8$   & $8/7$ \\
(m) & $2F_8$        & $1/2$ \\
\hline
\end{tabular} \qquad
\begin{tabular}[t]{l|l|l}
\multicolumn{3}{c}{$J=\{(3, 1), (5, 2)\}$} \\
\hline
\multicolumn{1}{c|}{case} &
\multicolumn{1}{c|}{$2E_S$} & \multicolumn{1}{c}{$(H |_S \cdot 2E |_S)$} \\
\hline
(a) & $2F_m$ & $4/m$ \\
(b) & $2F_4$ & $8/15$ \\
(c) & $2F_3$ & $1/3$ \\
(c) & $4F_3$ & $4/3$ \\
(d) & $2F_4$ & $1/5$ \\
(d) & $4F_4$ & $4/5$ \\
(e) & $2F_5$ & $2/15$ \\
(e) & $4F_5$ & $8/15$ \\
\hline
\end{tabular}
\end{center}
$(H |_S \cdot 2E |_S)=8E^3$ (\ref{eqn:key.value}),
which is $8/7$ (\ref{txt:typeI}(\ref{itm:typeI.7})),
$8/15$ (\ref{txt:typeI}(\ref{itm:typeI.15})) respectively.
Thus we have the following.
\end{Text}

\begin{Theorem}\label{thm:typeI}
If $f$ is of type I then $P$ is c$E_7$ or c$E_8$.
\end{Theorem}

\begin{Remark}\label{rem:typeI}
Actually $P$ is c$E_7$ when $J=\{(7, 3)\}$.
By (\ref{txt:list.typeI}) it suffices to exclude the case where
$P \in S_X$ is of type $E_8$ and $Z=2F_7+3F_8$.
In the case $S \cap E$ defines the $1$-cycle $F_7+F_8$ on $S$,
and $(E \cdot F_7)$, $(E \cdot F_8)$ must be $-1/7$.
On the other hand using the dual graph for the minimal resolution
(\ref{txt:graph}) we obtain
$(E |_S \cdot F_7)=(F_7+F_8 \cdot F_7)_S=-3/7$ and
$(E |_S \cdot F_8)=(F_7+F_8 \cdot F_8)_S=1/7$.
It is a contradiction.
\end{Remark}

\begin{Example}\label{exl:I}
There exist examples of type I. The weighted blowup of
$o \in (x_1^2+x_2^3+x_2x_3^3+x_4^7=0) \subset \bC^4$ with weights
$\wt(x_1,x_2,x_3,x_4)=(7,5,3,2)$
is an example of type I whose $P$ is c$E_7$,
and the weighted blowup of
$o \in (x_1^2+x_2^3+x_3^5+x_4^7=0) \subset \bC^4$ with weights
$\wt(x_1,x_2,x_3,x_4)=(7,5,3,2)$ is an example of type I whose $P$ is c$E_8$.
\end{Example}

\begin{Text}
On the other hand if $f$ is of type II or III,
then by (\ref{thm:el.irred}), (\ref{thm:II.III.elephant})
and (\ref{rem:precise}),
we obtain the dual graph for the partial resolution $S \to S_X$ below.
\begin{enumerate}
\item (type IIa or ${\textrm{IIb}}^{\vee}$) The dual graph is
\begin{align*}
\bullet - \ast
\begin{cases}
- \underset{1}{\circ} \\
\ \vdots \\
- \underset{1}{\circ},
\end{cases}
\end{align*}
where $\bullet$ denotes $H |_S=D$,
the attached numbers are the coefficients in $E |_S=Z$,
and $\ast$ means that a Du Val singularity of type
\begin{enumerate}
\item $A_{r-1}$ if $f$ is of type IIa,
\item $A_{2r-1}$ if $f$ is of type ${\textrm{IIb}}^{\vee}$ with
$J=\{(r, 1), (r, 1)\}$, $r \ge 5$,
\item $A_5$ or $E_6$
if $f$ is of type ${\textrm{IIb}}^{\vee}$ with $J=\{(3, 1), (3, 1)\}$,
\item $D_5$
if $f$ is of type ${\textrm{IIb}}^{\vee}$ with $J=\{(2, 1), (4, 1)\}$,
\end{enumerate}
appears there respectively.
The number of exceptional curves $\circ$ is at most $4$ in type IIa
and at most $a$ in type ${\textrm{IIb}}^{\vee}$.
This bound comes from the value of $(H \cdot [S \cap E])=aE^3$.
\item (type ${\textrm{IIb}}^{\vee \vee}$ or III)
Set $r_1=1$, $r_2=r$ if $f$ is of type III.
$P \in S_X$ is of type $A_n$ because $H |_S$ intersects $E |_S$
at two points.
The dual graph is
\begin{align*}
\bullet - \ast - \underset{1}{\circ} - \star - \bullet \quad \textrm{or} \quad
\bullet - \ast - \underset{1}{\circ} - \underset{1}{\circ} - \star - \bullet,
\end{align*}
where $\bullet$ denotes $H |_S=D$,
the attached numbers are the coefficients in $E |_S=Z$,
and $\ast$, $\star$ means that a Du Val singularity of type
$A_{r_1-1}$, $A_{r_2-1}$ appears there respectively.
Since all the components of $(S \cap E)_{\red}$ are numerically proportional,
the number of exceptional curves $\circ$ is one or two.

\end{enumerate}
The following list shows all the possibilities of the partial resolution
$S \to S_X$ up to permutation in terms of $F_i$ in (\ref{txt:graph}).

\begin{center}
\begin{tabular}{l|l|l|l}
\hline
\multicolumn{1}{c|}{type of $f$} & \multicolumn{1}{c|}{$S_X$} &
\multicolumn{1}{c|}{$Z$} & \multicolumn{1}{c}{$(H |_S \cdot E |_S)$} \\
\hline
IIa                     & $D_r$     & $F_r$                 & $4/r$ \\
IIa                     & $D_{r+1}$ & $F_r+F_{r+1}$         & $4/r$ \\
IIa                     & $D_{r+1}$ & $F_1+F_{r+1}$         & $1+1/r$ \\
IIa                     & $D_{r+2}$ & $F_1+F_{r+1}+F_{r+2}$ & $1+1/r$ \\
${\textrm{IIb}}^{\vee}$, $J=\{(r, 1), (r, 1)\}$ &
$D_{2r}$ & $F_{2r}$            & $2/r$ \\
${\textrm{IIb}}^{\vee}$, $J=\{(r, 1), (r, 1)\}$ &
$D_{2r+1}$ & $F_{2r}+F_{2r+1}$ & $2/r$ \\
${\textrm{IIb}}^{\vee}$, $J=\{(r, 1), (r, 1)\}$ &
$D_{2r+1}$ & $F_1+F_{2r+1}$    & $1+1/2r$ \\
${\textrm{IIb}}^{\vee}$, $J=\{(3, 1), (3, 1)\}$ &
$E_7$    & $F_6$    & $2/3$ \\
${\textrm{IIb}}^{\vee}$, $J=\{(2, 1), (4, 1)\}$ &
$D_6$    & $F_1$    & $1$ \\
${\textrm{IIb}}^{\vee}$, $J=\{(2, 1), (4, 1)\}$ &
$E_6$    & $F_2$    & $3/4$ \\
${\textrm{IIb}}^{\vee \vee}$ & $A_{r_1+r_2-1}$ & $F_{r_1}$ & $1/r_1+1/r_2$ \\
${\textrm{IIb}}^{\vee \vee}$ & $A_{r_1+r_2}$   &
$F_{r_1}+F_{r_1+1}$ & $1/r_1+1/r_2$ \\
III                     & $A_r$           & $F_r$         & $1+1/r$ \\
III                     & $A_{r+1}$       & $F_r+F_{r+1}$ & $1+1/r$ \\
\hline
\end{tabular}
\end{center}
$(H |_S \cdot E |_S)=aE^3$ (\ref{eqn:key.value}).
Combining it with (\ref{thm:num.classn}) we have the following theorem.
\end{Text}

\begin{Theorem}\label{thm:typeIIandIII}
Assume that $f$ is of type II or III.
Then possible types of Du Val singularities on $S$ and $S_X$ are as follows.
\begin{center}
\begin{tabular}{l|l|l}
\hline
\multicolumn{1}{c|}{type of $f$}   &
\multicolumn{1}{c|}{type of $S_X$} & \multicolumn{1}{c}{type of $S$} \\
\hline
IIa  & $D_r$ or $D_{r+1}$                       & $A_{r-1}$ \\
${\textrm{IIb}}^{\vee}$, $J=\{(r, 1), (r, 1)\}$
& $D_{2r}$ or $D_{2r+1}$ & $A_{2r-1}$ \\
${\textrm{IIb}}^{\vee}$, $J=\{(3, 1), (3, 1)\}$ & $E_7$    & $E_6$ \\
${\textrm{IIb}}^{\vee}$, $J=\{(2, 1), (4, 1)\}$ & $E_6$    & $D_5$ \\
${\textrm{IIb}}^{\vee \vee}$
& $A_{r_1+r_2-1}$ or $A_{r_1+r_2}$  & $A_{r_1-1}$ and $A_{r_2-1}$\\
III  & $A_r$ or $A_{r+1}$                       & $A_{r-1}$ \\
\hline
\end{tabular}
\end{center}
Furthermore the exceptional locus of the partial resolution $S \to S_X$
is irreducible unless the type of $S_X$ is that after
the word ``or" in the above table.
\end{Theorem}

\begin{Remark}\label{rem:typeIIandIII}
In types ${\textrm{IIb}}^{\vee \vee}$ and III, the type of $S_X$ has to be
in fact $A_{r_1+r_2-1}$, $A_r$ respectively by \cite{Ka01}, \cite{Ka01p} and
(\ref{thm:cAn}) proved in Section \ref{sec:cAn}.
\end{Remark}

\begin{Example}\label{exl:rr}
Let $f$ be the weighted blowup of
$o \in (x_1^2+x_2^2x_3+x_3^{2r}+x_4^r=0) \subset \bC^4$ with weights
$\wt(x_1,x_2,x_3,x_4)=(r,r,1,2)$, where $r \ge 3$ is an odd integer.
$P$ is c$D_{r+1}$ and
$f$ is of type ${\textrm{IIb}}^{\vee}$ with $J=\{(r, 1), (r, 1)\}$.
The exceptional locus of the partial resolution $S \to S_X$ is reducible.
\end{Example}

\begin{Example}\label{exl:24}
Let $f$ be the weighted blowup of
$o \in (x_1^2+x_2^2x_3+x_3^3+x_4^3=0) \subset \bC^4$ with weights
$\wt(x_1,x_2,x_3,x_4)=(3,1,4,2)$.
$P$ is c$D_4$ and
$f$ is of type ${\textrm{IIb}}^{\vee}$ with $J=\{(2, 1), (4, 1)\}$.
\end{Example}

\section{Divisorial contractions to c$\mathbf{A_n}$ points}\label{sec:cAn}
We begin with introducing a general method to determine $f$.

\begin{Lemma}\label{lem:gen.mtd}
Let $f \colon (Y \supset E) \to (X \ni P)$ be a germ of a
three dimensional divisorial contraction to a cDV point $P$
and let $K_Y=f^*K_X+aE$.
$X$ is identified with a hypersurface in $\bar{X}=\bC^4=
x_1x_2x_3x_4\textrm{-}\mathrm{space}$,
\begin{align*}
P \in X \cong o \in (\phi=0) \subset \bar{X}=\bC^4=
x_1x_2x_3x_4\textrm{-}\mathrm{space}.
\end{align*}
Let $m_i$ be the multiplicity of $\divis x_i$ along $E$, that is,
the largest integer such that $x_i \in f_*\cO_Y(-m_iE)$.
We assume that $m_1, m_2, m_3, m_4$ have no common factors.
Let $d$ be the weighted order of $\phi$ with respect to weights
$\wt(x_1, x_2, x_3, x_4)=(m_1, m_2, m_3, m_4)$ and decompose $\phi$ as
\begin{align*}
\phi=\phi_d(x_1, x_2, x_3, x_4)+\phi_{>d}(x_1, x_2, x_3, x_4).
\end{align*}
where $\phi_d$ be the weighted homogeneous part of weighted degree $d$
and $\phi_{>d}$ be the part of weighted degree greater than $d$.
Set $c=m_1+m_2+m_3+m_4-1-d$.
Let $\bar{g} \colon (\bar{Z} \supset \bar{F}) \to (\bar{X} \ni P)$
be the weighted blowup of $\bar{X}$ with weights
$\wt(x_1, x_2, x_3, x_4)=(m_1, m_2, m_3, m_4)$,
let $\bar{F}$ be its exceptional divisor, and
let $\bar{D_i}$ be the strict transform on $\bar{Z}$
of the hyperplane $x_i=0$ in $\bar{X}$.
Let $Z$ be the strict transform of $X$ on $\bar{Z}$,
and let $g \colon Z \to X$ be the induced morphism.
Assume that
\begin{enumerate}
\item $\bar{F} \cap Z$ defines an irreducible reduced $2$-cycle $F$ on $Z$,
and $F \not\subseteq \Sing Z$,\label{itm:2cycle}
\item $\dim(\Sing\bar{Z} \cap Z) \le 1$,\label{itm:dim.1}
\item $\Sing Z \subseteq \bigcup_{1 \le i \le 4} \bar{D_i}$,
and\label{itm:sing.cond}
\item
$F \not\subseteq \bigcup_{1 \le i \le 4} \bar{D_i}$.\label{itm:sing.cond2}
\end{enumerate}
Then $c \le a$, and the equality holds if and only if $f \cong g$ over $X$.
\end{Lemma}

\begin{proof}
$Z$ is $\mathrm{R}_1$ from (\ref{itm:2cycle}), and is Cohen-Macaulay
since $Z \subset \bar{Z}$ is locally a cyclic quotient of a hypersurface
in $\bC^4$. Thus $Z$ is normal.
From (\ref{itm:dim.1}) we can use the adjunction formula and
have $K_Z=(K_{\bar{Z}}+Z) |_Z=g^*K_X+cF$.
We note that $-F$ is $f$-ample and $\bQ$-Cartier.

Consider the centre on $Z$ of the valuation corresponding to $E$.
Because the multiplicity of $\divis x_i$ along $E$
equals to that along $F$,
this centre is not contained in $\bigcup_{1 \le i \le 4} \bar{D_i}$
by $\bigcap_{1 \le i \le 4} \bar{D_i}=\emptyset$ and (\ref{itm:sing.cond2}),
whence it intersects the smooth locus of $Z$ by (\ref{itm:sing.cond}).
Thus the discrepancy $a$ of $E$ with respect to $K_X$ is greater than
or equal to $c$, that of $F$, and the equality holds if and only if
$f \cong g$ over $X$ by \cite[Lemma 3.4]{Ka01}.
\end{proof}

\begin{Text}
Let $f \colon (Y \supset E) \to (X \ni P)$ be a germ of a divisorial
contraction to a c$A_n$ point $P$ ($n \ge 2$).
$X$ is identified with a hypersurface $\phi=0$ in $\bC^4=
x_1x_2x_3x_4\textrm{-}\mathrm{space}$ as follows.
\begin{align}\label{eqn:hyp.surf.cA}
P \in X \cong o \in (\phi=x_1x_2+g(x_3, x_4)=0) \subset \bC^4=
x_1x_2x_3x_4\textrm{-}\mathrm{space},
\end{align}
where the total order of $g$ with respect to $x_3, x_4$ is $n+1$.
From (\ref{rem:typeIV}) and (\ref{thm:typeI}),
$f$ is of type O, II or III.
Our goal in this section
is the classification of these contractions (\ref{thm:cAn}).
\end{Text}

\begin{Text}\label{txt:coord.change}
We start with an identification (\ref{eqn:hyp.surf.cA}).
Set $m_i$ as the multiplicity of $\divis x_i$ along $E$.
By coordinates change in which $x_4 \mapsto x_4+x_i$ ($i=1$, $2$ or $3$)
if necessary, we may assume that $m_4=1$.
By coordinates change $x_3 \mapsto x_3+h(x_4)$ if necessary
we may furthermore assume that $x_3+h(x_4) \not\in f_*\cO_Y(-(m_3+1)E)$
for any $h$. Decompose $\phi$ as
\begin{align*}
\phi=\phi_{\le m_1+m_2}(x_1, x_2, x_3, x_4)
+\phi_{>m_1+m_2}(x_1, x_2, x_3, x_4),
\end{align*}
where $\phi_{\le m_1+m_2}$ is the part of weighted degree less than or equal to
$m_1+m_2$ and $\phi_{> m_1+m_2}$ is the part of
weighted degree greater than $m_1+m_2$
with respect to weights $\wt(x_1, x_2, x_3, x_4)=(m_1, m_2, m_3, 1)$.

We focus on $\phi_{\le m_1+m_2}=x_1x_2+h(x_3, x_4)$. Of course $h \neq 0$.
The homogeneous part $h_0$ in $h$ of the minimal weighted degree $d_0$
with respect to weights
$\wt(x_3, x_4)=(m_3, 1)$ always decomposes into a product of
$x_4^{d_0-m_3\lfloor \frac{d_0}{m_3} \rfloor}$
and $\lfloor \frac{d_0}{m_3} \rfloor$ linear combinations of
$x_3$ and $x_4^{m_3}$. Hence $h_0 \not\in f_*\cO_Y(-(d_0+1)E)$
because of the assumption $x_3+h(x_4) \not\in f_*\cO_Y(-(m_3+1)E)$ for any $h$.
Therefore $d_0$ must be $m_1+m_2$.

We can see that this choice of coordinates satisfies the assumptions in
(\ref{lem:gen.mtd}), whence
\begin{align*}
a \ge m_1+m_2+m_3+1-1-(m_1+m_2)=m_3,
\end{align*}
and the equality holds if and only if $f \cong g$ in (\ref{lem:gen.mtd})
over $X$.
\end{Text}

In particular we say the following in type O.

\begin{Lemma}\label{lem:classnO}
If $f$ is of type O, then $f$ is the weighted blowup of $X$
with weights $\wt(x_1, x_2, x_3, x_4)=(i, n+1-i, 1, 1)$
\textup{(}$1 \le i \le n$\textup{)} in \textup{(\ref{eqn:hyp.surf.cA})}.
\end{Lemma}

\begin{Lemma}\label{lem:excludeIIa}
$f$ is not of type IIa with $a=4$.
\end{Lemma}

\begin{proof}
If $f$ is of type IIa with $a=4$, then
$J=\{(5, 2)\}$ from (\ref{thm:localIIa}), whence $Y$ has exactly one
non-Gorenstein point and it is a quotient singularity of type
$\frac{1}{5}(1, -1, 3)$.
Consider a birational morphism $h \colon Z \to Y$
such that $Z$ has an $h$-exceptional divisor $F$
whose discrepancy with respect to $K_X$ is $1$.
We can choose $n$ different $F$ as valuations corresponding to
exceptional divisors of the weighted blowups of $X$ with weights
$\wt(x_1, x_2, x_3, x_4)=(i, n+1-i, 1, 1)$ ($1 \le i \le n$)
in (\ref{eqn:hyp.surf.cA}).
Write
\begin{align*}
K_Z&=h^*K_Y+bF+(\textrm{others}),\\
h^*E&=E_Z+mF+(\textrm{others}),
\end{align*}
where $E_Z$ is the strict transform of $E$.
Then $K_Z=h^*(f^*K_X+4E)+bF+(\textrm{others})=h^*f^*K_X+4E_Z+
(b+4m)F+(\textrm{others})$, whence $b=m=1/5$ by $5b, 5m \in \bZ$.
Hence $F$ is determined uniquely because $Y$ has a quotient singularity of
index $5$ as the unique non-Gorenstein point.
It contradicts that we can choose $n$ different $F$.
\end{proof}

The case where $f$ is of type IIa with $a=2$, type IIb or III remains.
We can see the next lemma by following \cite[Claim 6.13]{Ka01p} faithfully.

\begin{Lemma}\label{lem:calculate}
Let $f \colon (Y \supset E) \to (X \ni P)$ be a germ of a
three dimensional divisorial contraction to a singular cDV point $P$,
and let $K_Y=f^*K_X+aE$.
Assume that $f$ is of type IIb or III, and
set $r_1=1$ if $f$ is of type III.
$P \in X$ is identified as
\begin{align*}
P \in X \cong o \in (\phi=0) \subset \bC^4=
x_1x_2x_3x_4\textrm{-}\mathrm{space}.
\end{align*}
Assume that the multiplicity of $\divis x_4$ along $E$ is $1$.
Then
\begin{enumerate}
\item for $1 \le j \le \min\{r_1, a\}$,
\begin{align*}
f_*\cO_Y(-jE)=(x_{1j}, x_{2j}, x_{3j}, x_4^i)\cO_X
\end{align*}
for some $x_{ij}=x_i+p_{ij}(x_4)$ \textup{(}$i=1, 2, 3$\textup{)}, where
\begin{align*}
p_{ij} \in \bigoplus_{k=1}^{j-1}\bC x_4^k \subset \bC[x_4].
\end{align*}
\item Assume that $r_1 < a$.
\begin{enumerate}
\item  The kernel of the linear map
\begin{multline*}
\bC x_{1r_1} \oplus \bC x_{2r_1} \oplus \bC x_{3r_1} \to\\
f_*\cO_Y(-r_1E)/(f_*\cO_Y(-(r_1+1)E)+(x_4^{r_1})\cO_X)
\end{multline*}
is of dimension $2$.
\item We assume that $x_{3r_1} \not\in f_*\cO_Y(-(r_1+1)E)+(x_4^{r_1})\cO_X$.
Under this situation, for $r_1 < j \le a$,
\begin{align*}
f_*\cO_Y(-jE)=(x_{1j}, x_{2j})\cO_X+
\sum_{(k,l) \in \bigcup_{s \ge j} J_s}(x_{3r_1}^kx_4^l)\cO_X
\end{align*}
for some $x_{ij}=x_i+p_{ij}(x_{3r_1}, x_4)$ \textup{(}$i=1, 2$\textup{)},
where
\begin{align*}
p_{ij} &\in \bigoplus_{(k,l) \in \bigcup_{1\le s<j}J_s}
\bC x_{3r_1}^kx_4^l \subset \bC[x_{3r_1}, x_4],\\
J_s &= \{(k,l) \in \bZ_{\ge 0}^2 \mid r_1k+l = s\}.
\end{align*}
\end{enumerate}
\end{enumerate}
\end{Lemma}

\begin{Remark}\label{rem:sim}
The similar result holds even if $f$ is of type IIa with $a=2$
from $d(-1)=\chi(\cQ_{-1})=1$ (\ref{rem:di}(\ref{itm:diIIa})).
Assume that $f$ is of type IIa with $a=2$ and
the multiplicity of $\divis x_4$ along $E$ is $1$.
Then $x_1, x_2, x_3 \in f_*\cO_Y(-2E)+(x_4^2)\cO_X$.
\end{Remark}

\begin{Text}
Assume that $f$ is of type IIa with $a=2$ or 
type ${\textrm{IIb}}^{\vee}$ with $J=\{(r, 1), (r, 1)\}$.
Then $m_3 \ge a$ in (\ref{txt:coord.change})
by (\ref{lem:calculate}), (\ref{rem:sim}).
Assume that $f$ is of type ${\textrm{IIb}}^{\vee \vee}$ or III,
and set $r_1=1, r_2=r$ if $f$ is of type III.
We start with an identification (\ref{eqn:hyp.surf.cA}).
Set $m_i$ as the multiplicity of $\divis x_i$ along $E$.
We may assume that $m_4=1$ by a coordinate change.
Because a general element in $f_*\cO_Y(-aE)$ gives a Du Val singularity
$P \in S_X$ at worst of type $A_{r_1+r_2}$ (\ref{thm:typeIIandIII}),
this element must be of form $cx_3+\cdots$ with $c \in \bC^\times$
by (\ref{lem:calculate}).
Hence by coordinates change such that
$x_3 \mapsto x_3+h(x_1, x_4)$ or $x_3+h(x_2, x_4)$
if necessary we may assume that $x_3 \in f_*\cO_Y(-aE)$ (\ref{lem:calculate}).
Therefore $m_3 \ge a$ in (\ref{txt:coord.change}).
Hence we have the following from (\ref{txt:coord.change}).
\end{Text}

\begin{Proposition}\label{prp:classnIIandIII}
Assume that $f$ is of type IIb or III.
If $f$ is not of type ${\textrm{IIb}}^{\vee}$ with $J=\{(2, 1), (4, 1)\}$,
then $f$ is the weighted blowup of $X$
with weights $\wt(x_1, x_2, x_3, x_4)=(m_1, m_2, a, 1)$ with respect
to the coordinates in \textup{(\ref{txt:coord.change})}.
\end{Proposition}

Finally we have to consider the exceptional case where
$f$ is of type ${\textrm{IIb}}^{\vee}$ with $J=\{(2, 1), (4, 1)\}$.

\begin{Lemma}\label{lem:excepIIb24}
If $f$ is of type ${\textrm{IIb}}^{\vee}$ with $J=\{(2, 1), (4, 1)\}$,
then
\begin{enumerate}
\item $P$ is a c$A_2$ point isomorphic to
$o \in (x_1x_2+x_3^3+g_{\ge4}(x_3, x_4)=0)
\subset \bC^4=x_1x_2x_3x_4\textrm{-}\mathrm{space}$,
where the total order of $g_{\ge4}$ with respect to $x_3, x_4$
is greater than or equal to $4$,\label{itm:cA2}
\item $Y$ has exactly one non-Gorenstein point $Q$, which is
isomorphic to $o \in (y_1^2+y_2^2+y_3^2+y_4^3=0)$ in the quotient space of
$\bC^4=y_1y_2y_3y_4\textrm{-}\mathrm{space}$ divided by $\bZ/(4)$ with
weights $\wt(y_1, y_2, y_3, y_4)=(1, 3, 3, 2)$, and\label{itm:descQ}
\item $K_Y=f^*K_X+3E$.\label{itm:a=3}
\end{enumerate}
\end{Lemma}

\begin{proof}
From (\ref{rem:precise}), $a=3$ (\ref{itm:a=3}) and
$Y$ has exactly one non-Gorenstein point $Q$.
The germ $Q \in Y$ is described as in (\ref{rem:precise}(\ref{itm:A/4})).
Because the Gorenstein index $r$ of $Y$ is $4$ and $a=3$,
for any valuation $v$ with discrepancy $1$ with respect to $K_X$
the centre of $v$ on $Y$ is $Q$ and $v$ has discrepancy $1/4$
with respect to $K_Y$, as in the proof of (\ref{lem:excludeIIa}).
Since there are $n$ such valuations
as we mentioned in the proof of (\ref{lem:excludeIIa}),
$Q \in Y$ must have an isomorphism in (\ref{itm:descQ})
by the description in (\ref{rem:precise}(\ref{itm:A/4})) and
\cite[Theorem 7.4]{Ha99}.
In this case there exist exactly two valuations centred at $Q$ with
discrepancy $1/4$ with respect to $K_Y$ \cite[Theorem 7.4]{Ha99}.
Therefore $P$ has to be c$A_2$.

We take an isomorphism $P \in X \cong o \in
(x_1x_2+g_3(x_3, x_4)+g_{\ge4}(x_3, x_4)=0) \subset \bC^4$,
where $g_3$ is the part of degree $3$
and $g_{\ge4}$ is the part of degree greater than or equal to $4$.
To see (\ref{itm:cA2}) we have to show that $g_3$ is cubic,
but this comes from the property that there exists a hyperplane section
of $P \in X$ which is Du Val of type $E_6$ at $P$ (\ref{thm:typeIIandIII}).
\end{proof}

\begin{Example}\label{exl:cA2}
There exists an example of (\ref{lem:excepIIb24}).
The weighted blowup of
$o \in (x_1^2+x_2^2+x_3^3+x_1x_4^2=0) \subset \bC^4$ with weights
$\wt(x_1,x_2,x_3,x_4)=(4,3,2,1)$ is such an example.
\end{Example}

We have (\ref{lem:classnO}), (\ref{lem:excludeIIa}), (\ref{prp:classnIIandIII})
and (\ref{lem:excepIIb24}).
We can easily check the sufficient condition for $f$ to be a
divisorial contraction, and obtain the classification (\ref{thm:cAn}).

\small{\textsc{Graduate School of Mathematical Sciences,
the University of Tokyo, 3-8-1 Komaba, Meguro, Tokyo 153-8914, Japan}}\\
\textit{email address:}
\texttt{kawakita@ms.u-tokyo.ac.jp}
\end{document}